\documentclass[12pt]{article}
\usepackage{amsmath,amssymb,amsfonts,graphicx}

\newtheorem{theorem}{Theorem}
\newtheorem{state}{Statement}
\newtheorem{lemma}{Lemma}

\begin{document}


\begin{center}
{\Large\bf On the convergence of probabilities of the random
graphs' properties expressed by first-order formulae with a
bounded quantifier depth}\footnote{This work was carried out with
the support of the Russian Foundation for Basic Research (grant
12-01-00683)} \vspace{0.5cm}

M.E. Zhukovskii
\end{center}

\section{Introduction}

An asymptotic behavior of the probabilities of
Erd\H{o}s--R\'{e}nyi random graph first-order properties is
studied in the article. In this section we briefly describe the
history of the problem and introduce necessary definitions. At the
end of the section we formulate our main result.

Let $N\in\mathbb{N}$, $0\leq p\leq 1.$ Denote the set of all
undirected graphs without loops and multiple edges with a set of
vertices $V_N=\{1,...,N\}$ by $\Omega_N=\{G=(V_N,E)\}$. \emph{The
Erd\H{o}s--R\'{e}nyi random graph} (see \cite{Erdos}--\cite{Alon})
is a random element $G(N,p)$ of $\Omega_N$ with a distribution
${\sf P}_{N,p}$ on $\mathcal{F}_N=2^{\Omega_N}$ defined as
follows:
$$
 {\sf P}_{N,p}(G)=p^{|E|}(1-p)^{C_N^2-|E|}.
$$

{ \it The random graph obeys zero-one law with a class of
properties $\mathcal{C}$} if for any property $C\in\mathcal{C}$
either $\lim\limits_{N\rightarrow\infty}{\sf P}_{N,p}(C)=0$ or
$\lim\limits_{N\rightarrow\infty}{\sf P}_{N,p}(C)=1$.

The class of first-order properties is the most studied class in
this area. Such properties are expressed by first-order formulae
(see \cite{Veresh}, \cite{Usp}). These formulae are built of
predicate symbols $\sim,=$, logical connectivities
$\neg,\Rightarrow,\Leftrightarrow,\vee,\wedge$, variables
$x,y,x_1...$ and quantifiers $\forall,\exists$. Symbols
$x,y,x_1...$ express vertices of a graph. The relation symbol
$\sim$ expresses the property of two vertices to be adjacent. The
symbol $=$ expresses the property of two vertices being
coincident. We denote by $\mathcal{P}$ a class of functions
$p=p(N)$ such that the random graph $G(N,p)$ obeys zero-one law
with the class
$\mathcal{L}$ of all first-order properties. 
In 1969 by Y.V. Glebskii, D.I. Kogan, M.I. Liagonkii and V.A.
Talanov in \cite{Kogan} (and independently in 1976 R.Fagin in
\cite{Fagin}) proved that if
$$
 \forall \alpha>0 \,\,\,\, N^{\alpha}\min\{p,1-p\}\rightarrow\infty,
 \,\, N\rightarrow\infty,
$$
then $p\in\mathcal{P}$. Moreover in 1988 S. Shelah and J.H.
Spencer (see \cite{Shelah})  expanded the class of functions
$p(N)$ ``that follow the zero-one law''. They proved that the
functions $p=N^{-\alpha}$,
$\alpha\in\mathbb{R}\setminus\mathbb{Q}$, $\alpha\in(0,1)$, are in
$\mathcal{P}$. Surely $p=1-N^{-\alpha}\in\mathcal{P}$ when
$\alpha\in\mathbb{R}\setminus\mathbb{Q}$, $\alpha\in(0,1)$.

If $\alpha$ is rational, $0<\alpha\leq 1$ and $p=N^{-\alpha}$ then
$G(N,p)$ does not obey zero-one law (see \cite{Alon}).

Denote by $\mathcal{L}^{\infty}$,
$\mathcal{L}^{\infty}\supset\mathcal{L}$, a class of all
properties expressed by formulae containing infinite number of
conjunctions and disjunctions. A class $\mathcal{L}^{\infty}_k$,
$\mathcal{L}^{\infty}_k\subset\mathcal{L}^{\infty}$, containing
all properties expressed by formulae with quantifier depths
bounded by the number $k$, in the frame of zero-one laws was
considered by M.~McArtur in 1997 (see \cite{McArtur}). M.~McArthur
obtained zero-one laws with the class $\mathcal{L}^{\infty}_k$ for
the random graph $G(N,N^{-\alpha})$ with some rational $\alpha$
from $(0,1]$.

Finally, the random graph $G(N,p)$ does not obey zero-one law with
the class $\mathcal{L}$ if $p=N^{-\alpha}$ and $\alpha$ is
rational, $\alpha\in(0,1]$. At the same time the random graph
$G(N,p)$ obeys zero-one law with the class
$\mathcal{L}^{\infty}_k$ for some rational $\alpha\in(0,1]$.
Therefore it seems natural to consider the class
$\mathcal{L}_k=\mathcal{L}\cap\mathcal{L}_k^{\infty}$. In 2010
(see \cite{Zhuk2}, \cite{zhuk_dan}) we proved that if $k\geq 3$,
$\alpha\in(0,1/(k-2))$ the random graph $G(N,N^{-\alpha})$ obeys
zero-one law with the class $\mathcal{L}_k$. We also proved that
when $\alpha=1/(k-2)$ the random graph $G(N,N^{-\alpha})$ does not
obey zero-one law with this class. This result led us to the
following question. Do probabilities ${\sf
P}_{N,N^{-1/(k-2)}}$ of all properties from $\mathcal{L}_k$ converge?\\

Let us state the main result of the article.

\begin{theorem}
Let $k\geq 3$, $p=N^{-\alpha}$, $\alpha=\frac{1}{k-2}$. For any
property $L\in\mathcal{L}_k$ there exists
$\lim_{N\rightarrow\infty}{\sf P}_{N,p}(L)$. \label{main}
\end{theorem}

Here we prove Theorem \ref{main} for the case $k\geq 4$ only. The
case $k=3$ is much easier and its correctness can be proved by
using the same arguments as in Lemma \ref{classes} (see Section 5
and Subsection 8.4).

We give a proof of Theorem \ref{main} in Section 7. This proof is
based on a number of statements from Section 5 and Section 6. The
mentioned statements are proved in Section 8. The main statement
is Lemma 1 which is related to the Ehrenfeucht game (see Section
4). It plays a key role in proofs of zero-one laws and other
statements on first-order properties. We introduce all necessary
constructions in Section 3 which is divided into 4 subsections. We
describe its structure in the end of Section 2 which is devoted to
some important and well-known theorems on extensions of small
subgraphs in the random graph.

\section{Distribution of small subgraphs}

For an arbitrary graph $G$ denote by $v(G)$ and $e(G)$ the number
of its vertices and the number of its edges respectively. The
number $\rho(G)=\frac{e(G)}{v(G)}$ is called \emph{density of
$G$.} The graph $G$ is called \emph{balanced} if for any subgraph
$H\subseteq G$ the inequality $\rho(H)\leq \rho(G)$ holds. The
graph $G$ is \emph{strictly balanced} if for any subgraph
$H\subset G$ the inequality $\rho(H)<\rho(G)$ holds.

Let us describe a problem studied by J.H.~Spencer in 1990 (see
\cite{Alon}, \cite{Spencer}). Consider graphs
$H,G,\widetilde{H},\widetilde{G}$. Let $V(H)=\{x_1,...,x_k\}$,
$V(G)=\{x_1,...,x_l\}$,
$V(\widetilde{H})=\{\widetilde{x}_1,...,\widetilde{x}_k\}$,
$V(\widetilde{G})=\{\widetilde{x}_1,...,\widetilde{x}_l\}$,
$H\subset G$, $\widetilde{H}\subset\widetilde{G}$ (therefore,
$k<l$). The graph $\widetilde{G}$ is called a
$(G,H)$-\emph{extension of the graph $\widetilde{H}$} if
$$
 \{x_{i_1},x_{i_2}\}\in E(G)\setminus E(H) \Rightarrow
 \{\widetilde{x}_{i_1},\widetilde{x}_{i_2}\}\in
 E(\widetilde{G})\setminus E(\widetilde{H}).
$$
If
$$
 \{x_{i_1},x_{i_1}\}\in
 E(G)\setminus E(H) \Leftrightarrow
 \{\widetilde{x}_{i_1},\widetilde{x}_{i_2}\}\in
 E(\widetilde{G})\setminus E(\widetilde{H})
$$
then we call $\widetilde{G}$ a \emph{strict extension}. Set
$$
 v(G,H)=|V(G)\setminus V(H)|, \,\,
 e(G,H)=|E(G)\setminus E(H)|,
$$
$$
 f_{\alpha}(G,H)=v(G,H)-\alpha e(G,H).
$$
Fix an arbitrary $\alpha>0$. If the inequality $f(S,H)>0$ holds
for any graph $S$ such that $H\subset S\subseteq G$ then the pair
$(G,H)$ is called \emph{$\alpha$-safe} (see \cite{Janson},
\cite{Alon}, \cite{Spencer}). If the inequality $f(G,S)<0$ holds
for any graph $S$ such that $H\subseteq S\subset G$ then the pair
$(G,H)$ is called {\it $\alpha$-rigid} (see \cite{Janson},
\cite{Alon}). The pair $(G,H)$ is called {\it $\alpha$-neutral} if
the following three properties hold. For any vertex $x$ of the
graph $H$ there exists a vertex of $V(G)\setminus V(H)$ adjacent
to $x$; $f_{\alpha}(S,H)>0$ for any graph $S$ such that $H\subset
S\subset G$; $f_{\alpha}(G,H)=0$.

Introduce a definition of a maximal pair. Let
$\widetilde{H}\subset\widetilde{G}\subset\Gamma$, $T\subset K$,
$|V(T)|\leq|V(\widetilde{G})|.$ The pair
$(\widetilde{G},\widetilde{H})$ is called \emph{$(K,T)$-maximal in
$\Gamma$} if for any subgraph $\widetilde{T}$ of $\widetilde{G}$
such that $|V(\widetilde{T})|=|V(T)|$ and
$\widetilde{T}\cap\widetilde{H}\neq\widetilde{T}$ the following
property holds. There is no $(K,T)$-extension $\widetilde{K}$ of
$\widetilde{T}$ in
$\Gamma\setminus(\widetilde{G}\setminus\widetilde{T})$ such that
each vertex of $V(\widetilde{K})\setminus V(\widetilde{T})$ is not
adjacent to any vertex of $V(\widetilde{G})\setminus
V(\widetilde{T})$.

Let $\alpha\in(0,1]$. Let a pair $(G,H)$ be $\alpha$-safe. Let
$V(H)=\{x_1,...,x_k\}$, $V(G)=\{x_1,...,x_l\}$. Denote  a set of
all $\alpha$-rigid pairs $(K_i,T_i)$ such that $|V(T_i)|\leq
|V(G)|,$ $|V(K_i)\setminus V(T_i)|\leq r$ by
$\Sigma^{\mbox{\small{rigid}}}(r)$. Consider a set
$\Sigma^{\mbox{\small{neutral}}}(r)$ of all $\alpha$-neutral pairs
$(K_i,T_i)$ such that $|V(T_i)|\leq |V(G)|,$ $|V(K_i)\setminus
V(T_i)|\leq r.$\\

Consider the random graph $G(N,p)$. Let $H\subset G$,
$V(H)=\{x_1,...,x_k\}$, $V(G)=\{x_1,...,x_l\}$,
$\widetilde{x}_1,...,\widetilde{x}_k\in V_N$. Define a random
variable $N_{(G,H)}(\widetilde{x}_1,...,\widetilde{x}_k)$ on the
probability space $(\Omega_N,\mathcal{F}_N,{\sf P}_{N,p})$ as
follows. The random variable assigns a number of all
$(G,H)$-extensions induced on the set
$\{\widetilde{x}_1,...,\widetilde{x}_k\}$ in $\mathcal{G}$ to a
graph $\mathcal{G}$ from $\Omega_N$. A graph $X$ is called \emph{a
subgraph of a graph $Y$ induced on a set $S\subset V(Y)$} if
$V(X)=S$ and for any vertices $x,y\in S$ the property  $\{x,y\}\in
E(X)\Leftrightarrow\{x,y\}\in E(Y)$ holds. Let us give a formal
definition of $N_{(G,H)}(\widetilde{x}_1,...,\widetilde{x}_k)$.
Let $W\subset V_N
\setminus\{\widetilde{x}_1,...,\widetilde{x}_k\},$  $|W|=l-k.$ If
there is a numeration of elements of the set $W$ by numbers
$k+1,k+2,...,l$ such that the graph
$\mathcal{G}|_{\{\widetilde{x}_1,...,\widetilde{x}_l\}}$ is a
$(G,H)$-extension of a graph
$\mathcal{G}|_{\{\widetilde{x}_1,...,\widetilde{x}_k\}}$ then we
set $I_W(\mathcal{G})=1.$ Otherwise we set $I_W(\mathcal{G})=0.$
The random variable
$N_{(G,H)}(\widetilde{x}_1,...,\widetilde{x}_k)$ is defined by the
equality
$$
N_{(G,H)}(\widetilde{x}_1,...,\widetilde{x}_k)=
\sum\limits_{W\subset V_N
\setminus\{\widetilde{x}_1,...,\widetilde{x}_k\},\, |W|=l-k}I_W.
$$

\begin{theorem} [\cite{Spencer}]
Let $p=N^{-\alpha}$. Let a pair $(G,H)$ be $\alpha$-safe. Then
$$
\lim\limits_{N\rightarrow\infty}{\sf P}_{N,p}(\forall
\widetilde{x}_1,...,\widetilde{x}_k \,\,
\left|N_{(G,H)}(\widetilde{x}_1,...,\widetilde{x}_k)- {\sf
E}_{N,p}N_{(G,H)}(\widetilde{x}_1,...,\widetilde{x}_k)\right|\leq
\varepsilon {\sf
E}_{N,p}N_{(G,H)}(\widetilde{x}_1,...,\widetilde{x}_k))=1
$$
for any $\varepsilon>0$. Here ${\sf E}_{N,p}$ is the expectation.
Moreover, ${\sf
E}_{N,p}N_{(G,H)}(\widetilde{x}_1,...,\widetilde{x}_k))=\Theta(N^{f(G,H)}).$
\label{safe_extensions}
\end{theorem}

In fact, the statement of this theorem means that almost surely
for any vertices $\widetilde{x}_1,...,\widetilde{x}_k$ the
relation
$$
 N_{(G,H)}(\widetilde{x}_1,...,\widetilde{x}_k)\sim
 {\sf E}_{N,p} N_{(G,H)}(\widetilde{x}_1,...,\widetilde{x}_k)
$$
holds. In such cases we will use this notation.

In addition to Theorem \ref{safe_extensions} J.H.~Spencer and
S.~Shelah (see \cite{Alon}, \cite{Shelah}) proved a result on a
number of maximal extensions of subgraphs in random graphs (in the
case of ``prohibited'' rigid pairs). In 2010 we extended this
result by considering ``prohibited'' neutral pairs (see
\cite{zhuk_extensions}).

Let us define new random variables and formulate the corresponding
results. Consider a random variable
$\widehat{N}^{\mbox{\small{rigid}}}_{(G,H),r}(\widetilde{x}_1,...,\widetilde{x}_k)$
such that if $\mathcal{G}\in\Omega_N$ then
$\widehat{N}^{\mbox{\small{rigid}}}_{(G,H),r}(\widetilde{x}_1,...,\widetilde{x}_k)[\mathcal{G}]$
is the number of strict $(G,H)$-extensions $\widetilde{G}$ of the
graph
$\widetilde{H}=\mathcal{G}|_{\{\widetilde{x}_1,...,\widetilde{x}_k\}}$
with the following property. For each pair
$(K_i,T_i)\in\Sigma^{\mbox{\small{rigid}}}(r)$ the pair
$(\widetilde{G},\widetilde{H})$ is $(K_i,T_i)$-maximal in
$\mathcal{G}$. First we formulate a result proved by
J.H.~Spencer and S.~Shelah in \cite{Shelah}.\\

\begin{theorem} [\cite{Shelah}]
Almost surely for any vertices
$\widetilde{x}_1,...,\widetilde{x}_k$
$$
 \widehat{N}^{\mbox{\small{rigid}}}_{(G,H),r}(\widetilde{x}_1,...,\widetilde{x}_k)\sim
 N_{(G,H)}(\widetilde{x}_1,...,\widetilde{x}_k) \sim
 {\sf
 E}_{N,p}
 \widehat{N}^{\mbox{\small{rigid}}}_{(G,H),r}(\widetilde{x}_1,...,\widetilde{x}_k)
 =\Theta\left(N^{f(G,H)}\right).
$$
\label{maximal_extensions_rigid}
\end{theorem}

Recall a result from \cite{zhuk_extensions}. Consider a random
variable
$\widehat{N}^{\mbox{\small{neutral}}}_{(G,H),r}(\widetilde{x}_1,...,\widetilde{x}_k)$
such that if $\mathcal{G}\in\Omega_N$ then
$\widehat{N}^{\mbox{\small{neutral}}}_{(G,H),r}(\widetilde{x}_1,...,\widetilde{x}_k)[\mathcal{G}]$
is the number of strict $(G,H)$-extensions $\widetilde{G}$ of the
graph
$\widetilde{H}=\mathcal{G}|_{\{\widetilde{x}_1,...,\widetilde{x}_k\}}$
with the following property. The pair
$(\widetilde{G},\widetilde{H})$ is $(K_i,T_i)$-maximal in
$\mathcal{G}$ for any
$(K_i,T_i)\in\Sigma^{\mbox{\small{neutral}}}(r)$.\\

\begin{theorem} [\cite{zhuk_extensions}]
Almost surely for any vertices
$\widetilde{x}_1,...,\widetilde{x}_k$
$$
 \widehat{N}^{\mbox{\small{neutral}}}_{(G,H),r}(\widetilde{x}_1,...,\widetilde{x}_k)\sim
 {\sf
 E}_{N,p}
 \widehat{N}^{\mbox{\small{neutral}}}_{(G,H),r}(\widetilde{x}_1,...,\widetilde{x}_k)
 =\Theta\left(N^{f(G,H)}\right).
$$
\label{maximal_extensions_neutral}
\end{theorem}

Let us proceed on to the proof of Theorem \ref{main}. An idea of
the proof is in the analysis of the probability of the existence
of a wining strategy for the second player called Duplicator in
the Ehrenfeucht game (see Section 4). In Section 3, all the
constructions which are necessary for the proof will be presented.

The next section consists of 4 subsections. The main constructions
used in the proof of Lemma 1 are introduced in Subsections 3.3,
3.4.
These constructions are maximal for all $\alpha$-neutral and
$\alpha$-rigid pairs that the first player called Spoiler can
build during $k$ rounds. In Subsection 3.2 the notion of a closure
$[A]_{\widehat{G}}$ for a subgraph $A$ of some graph $\widehat{G}$
is introduced.

This closure ``contains'' all $\alpha$-neutral pairs. So, in
Subsections 3.3, 3.4 graphs containing a maximal number of
$\alpha$-rigid pairs are constructed. Then closures of such graphs
are considered. In Subsection 3.1 all necessary pairs of graphs
are defined.

\section{Constructions}

Let $k\geq 4$ be natural. In what follows we assume
$\alpha=1/(k-2)$. So, we will write $f(G,H)$ instead of
$f_{\alpha}(G,H)$ everywhere below.

\subsection{Additional graphs}

Consider graphs
$H_1,H_2,G_1,G_2,G_3^{i_1,...,i_t},G_4,G_4^1,G_4^2$, where
$t\in\{1,...,k-2\}$, $i_1\in\{1,...,k-2\}$,
$i_2\in\{1,...,k-2\}\setminus\{i_1\}$, ...,
$i_t\in\{1,...,k-2\}\setminus\{i_1,...,i_{t-1}\}$.

In the following subsections we will use pairs of these graphs. In
fact, we are interested in the pairs $(G_1,H_1)$, $(G_2,H_2)$,
$(G_3^{i_1,...,i_t},H_2)$, $(G_4^1,H_1)$, $(G_4^2,H_1)$.

\begin{itemize}

\item[1)] The graph $G_1$ is complete, $V(G_1)=\{x_1,...,x_k\}$,
$H_1$ is an arbitrary graph on the set of vertices
$V(H_1)=\{x_1,...,x_{k-3}\}$.

\item[2)] Let
$$
 V(H_2)=\{x_1,...,x_{k-3},x_{k-2}\}, \,\,
 E(H_2) - \mbox{arbitrary};
$$
$$
V(G_2)=V(H_2)\cup
 \{x_{k-1}\}, \,\,
 E(G_2)= E(H_2)\cup\{\{x_1,x_{k-1}\},...,\{x_{k-2},x_{k-1}\}\}.
$$

\item[3)] Let $t\in\{1,...,k-2\}$, $i_1\in\{1,...,k-2\}$,
$i_2\in\{1,...,k-2\}\setminus\{i_1\}$, ...,
$i_t\in\{1,...,k-2\}\setminus\{i_1,...,i_{t-1}\}$. Consider graphs
$G_3^{i_1,...,i_t}$ defined by induction:
$$
 V(G_3^{i_1})=V(G_2)\cup\{x_{k}^{i_1}\},
$$
$$
 E(G_3^{i_1})=E(G_2)\cup
 \{\{x_1,x_k^{i_1}\},...,\{x_{k-1},x_k^i\}\}\setminus\{\{x_{i_1},x_k^{i_1}\}\};
$$
$$
 V(G_3^{i_1,...,i_t})=V(G_3^{i_1,...,i_{t-1}})\cup\{x_{k}^{i_1,...,i_t}\},
$$
$$
 E(G_3^{i_1,...,i_t})=E(G_3^{i_1,...,i_{t-1}})\cup
 \{\{x_1,x_{k}^{i_1,...,i_t}\},...,\{x_{k-1},x_{k}^{i_1,...,i_t}\}\}\setminus\{\{x_{i_t},x_{k}^{i_1,...,i_t}\}\}.
$$

\item[4)] Let
$$
 V(G_4)=V(H_1)\cup\{x_{k+1},x_{k+2},x_{k+3}\},
$$
$$
 V(G_4^1)=V(G_4)\cup\{x_{k+4}^1\},\,\,
 V(G_4^2)=V(G_4)\cup\{x_{k+4}^2,x_{k+5}^2\};
$$
$$
 E(G_4)=E(H_1)\cup\{\{x_1,x_{k+1}\},...,\{x_{k-4},x_{k+1}\},
 \{x_1,x_{k+2}\},...,\{x_{k-3},x_{k+2}\},
$$
$$
 \{x_1,x_{k+3}\},...,\{x_{k-3},x_{k+3}\},
 \{x_{k+1},x_{k+2}\},\{x_{k+1},x_{k+3}\},\{x_{k+2},x_{k+3}\}\};
$$
$$
 E(G_4^1)=E(G_4)\cup\{\{x_1,x_{k+4}^1\},...,\{x_{k-3},x_{k+4}^1\},
 \{x_{k+1},x_{k+4}^1\},\{x_{k+3},x_{k+4}^1\}\};
$$
$$
 E(G_4^2)=E(G_4)\cup\{\{x_1,x_{k+4}^2\},...,\{x_{k-3},x_{k+4}^2\},
 \{x_1,x_{k+5}^2\},...,\{x_{k-3},x_{k+5}^2\},
$$
$$
 \{x_{k+1},x_{k+4}^2\},\{x_{k+1},x_{k+5}^2\},
 \{x_{k+4}^2,x_{k+5}^2\}\}.
$$

\end{itemize}

Let $t\in\{1,...,k-2\}$, $i_1\in\{1,...,k-2\}$,
$i_2\in\{1,...,k-2\}\setminus\{i_1\}$, ...,
$i_t\in\{1,...,k-2\}\setminus\{i_1,...,i_{t-1}\}$. Consider the
set $S^{i_1,...,i_t}$ of all unordered collections of $k-2$
vertices from $V(G_3^{i_1,...,i_t})$. For each $U\subset
S^{i_1,...,i_t}$ consider the union of the graph
$G_3^{i_1,...,i_t}$ and all the $(G_2,H_2)$-extensions of graphs
$G_3^{i_1,...,i_t}|_{\mathbf{u}}$, $\mathbf{u}\in U$. Denote this
union by $G_3^{i_1,...,i_t}(U)$. Note that a union of a graph on
$k-2$ vertices with its $(G_2,H_2)$-extension is obtained by
adding one vertex adjacent to all its vertices. Let
$\mathcal{U}_i^{i_1,...,i_t}$ be a set of all subsets of
$S^{i_1,...,i_t}$ with the cardinality $i$.\\

Let us construct the closure of a graph.

\subsection{The closure $[A]_{\widehat{G}}$ in $\widehat{G}$ of a graph $A$}

Consider arbitrary vertices $\widehat{x}_1,...,\widehat{x}_{k-3}$.
Let $\widehat{G}$ be a graph. Let
$\widehat{x}_1,...,\widehat{x}_{k-3}$ be the vertices of the graph
$\widehat{G}$. Consider any graph $A\subset\widehat{G}$ on a set
of vertices
$V(A)=\{a_1,...,a_d\}$, $d\geq k-2$.\\

First of all let us note that for the graph $A$ there exist
several closures in the graph $\widehat{G}$. All these closures
are isomorphic.\\

We construct the graph $[A]_{\widehat{G}}$ in $k-1$ steps. Let $S$
be a set of all different unordered collections of $k-2$ vertices
from $V(A)$. Set $[A]_{\widehat{G}}= A$,
$\widehat{G}_1=\widehat{G}$.\\

{\bf The first step} is divided into $|S^{1,...,k-2}|$ parts.
Consider the first part of the step. Let
$\{a_{i_1},...,a_{i_{k-2}}\}\in S$,
$U\in\mathcal{U}^{1,...,k-2}_{|S^{1,...,k-2}|}=\{S^{1,...,k-2}\}$.
Assume that there exists an $(G_3^{1,...,k-2}(U),H_2)$-extension
$\widehat{Q}$ of $A|_{\{a_{i_1},...,a_{i_{k-2}}\}}$ in
$\widehat{G}_1$. We add only one such extension to the graph
$[A]_{\widehat{G}}$ and for all these extensions we remove
extenders from the graph $\widehat{G}_1$ (if a graph $X$ is an
extension of a graph $Y$ then we say that graph $X\setminus Y$ is
{\it an extender}).

Let the first $s$ parts, $s\leq |S^{1,...,k-2}|-1$, of the first
step of the graphs $[A]_{\widehat{G}},\widehat{G}_1$ construction
be done. Let us describe the $s+1$-th part. Let
$\{a_{i_1},...,a_{i_{k-2}}\}\in S$,
$U\in\mathcal{U}^{1,...,k-2}_{|S^{1,...,k-2}|-s}$. Assume that an
$(G_3^{1,...,k-2}(U),H_2)$-extension $\widehat{Q}$ of
$A|_{\{a_{i_1},...,a_{i_{k-2}}\}}$ in $\widehat{G}_1$ exists. We
add only one such extension to the graph $[A]_{\widehat{G}}$ and
for all these extensions we remove the extenders from the graph
$\widehat{G}_1$.\\

Let the $i$-th step, $i\leq k-3$, be done. Describe {\bf the
$i+1$-th step}. We divide this step into $|S^{1,...,k-2-i}|$
parts.

Let $\{a_{i_1},...,a_{i_{k-2}}\}\in S$,
$i_1,...,i_{k-2-i}\in\{1,...,k-2\}$ be an ordered collection of
different numbers, $U$ be a subset of
$\mathcal{U}^{i_1,...,i_{k-2-i}}_{|S^{i_1,...,i_{k-2-i}}|}$. Let
an $(G_3^{i_1,...,i_{k-2-i}}(U),H_2)$-extension $\widehat{Q}$ of
$A|_{\{a_{i_1},...,a_{i_{k-2}}\}}$ in $\widehat{G}_1$ exist. We
add only one such extension to the graph $[A]_{\widehat{G}}$ and
for all these extensions we remove the extenders from the graph
$\widehat{G}_1$.

Let the first $s$ parts, $s\leq |S^{i_1,...,i_{k-2-i}}|-1$, of the
$i+1$-th step be done, $\{a_{i_1},...,a_{i_{k-2}}\}\in S$. Let
$i_1,...,i_{k-2-i}\in\{1,...,k-2\}$ be an unordered collection of
different numbers. Let
$U\in\mathcal{U}^{i_1,...,i_{k-2-i}}_{|S^{i_1,...,i_{k-2-i}}|-s}$.
Assume that an $(G_3^{i_1,...,i_{k-2-i}}(U),H_2)$-extension
$\widehat{Q}$ of $A|_{\{a_{i_1},...,a_{i_{k-2}}\}}$ in
$\widehat{G}_1$ exists. We add only one such extension to the
graph $[A]_{\widehat{G}}$ and for all these extensions we remove
the extenders from the graph
$\widehat{G}_1$.\\

Describe the final {\bf $k-1$-th step}. Let
$\{a_{i_1},...,a_{i_{k-2}}\}\in S$. Assume that there exists a
$(G_2,H_2)$-extension $\widehat{Q}$ of
$A|_{\{a_{i_1},...,a_{i_{k-2}}\}}$ in $\widehat{G}_1$.
We add only one such extension to the graph $[A]_{\widehat{G}}$. The graph $[A]_{\widehat{G}}$ is constructed.\\

In the following two subsections we construct graphs
$X^l_{\widehat{G}}(\widehat{x}_1,...,\widehat{x}_{k-3})$,
$\widehat{X}^l_{\widehat{G}}(\widehat{x}_1,...,\widehat{x}_{k-3})$,
$X^l_j(\widehat{x}_1)$, $\widehat{X}^l_j(\widehat{x}_1)$, where
$l\in\{1,2,3,4,5\}$, $j$ is from some set $J$. The graphs
$\widehat{X}^l_j(\widehat{x}_1)$ are subgraphs of some graphs
$\widehat{G}_{i_j}$ which are chosen in such a way that these
subgraphs are ``different'' in some sense. The graphs
$\widehat{X}^l_{\widehat{G}}(\widehat{x}_1,...,\widehat{x}_{k-3})$
are subgraphs of the graph $\widehat{G}$ and built as $l$ grows
from $1$ to $5$. The graph $\widehat{X}^l_j(\widehat{x}_1)$ is the
union of the graphs
$\widehat{X}^l_{\widehat{G}_{i_j}}(\widehat{x}_1,\widehat{x}_2^i,...,\widehat{x}_{k-3}^i)$
over some sets of vertices
$\widehat{x}_{2}^i,...,\widehat{x}_{k-3}^i$ of the graph
$\widehat{G}_{i_j}$. Finally, graphs
$X^l_{\widehat{G}}(\widehat{x}_1,...,\widehat{x}_{k-3})$,
$X^l_j(\widehat{x}_1)$ are the unions of the closures of some
subgraphs of
$\widehat{X}^l_{\widehat{G}}(\widehat{x}_1,...,\widehat{x}_{k-3})$,
$\widehat{X}^l_j(\widehat{x}_1)$ respectively.

\subsection {Graphs
$X^l_{\widehat{G}}(\widehat{x}_1,...,\widehat{x}_{k-3})$,
$\widehat{X}^l_{\widehat{G}}(\widehat{x}_1,...,\widehat{x}_{k-3})$,
$l\in\{1,2,3,4,5\}$}

Assume that the graph $\widehat{G}$ considered in the previous
subsection does not contain subgraphs $W$ with $v(W)<k^3$ and
$\rho(W)>k-2$. We consider this restriction because of the
following reasonings. As it is mentioned in Section 2 the
constructions we build should contain a maximal number of rigid
pairs. Without a restriction on the density of subgraphs in
$\widehat{G}$ the number of such pairs in $\widehat{G}$ can be
arbitrarily large. We choose the number $k-2$ as in the random
graph $G(N,p)$ there are no subgraphs $W$ with
$v(W)$ bounded by a fixed number and $\rho(W)>k-2$.\\

We construct the subgraph
$\widehat{X}^1_{\widehat{G}}(\widehat{x}_1,...,\widehat{x}_{k-3})$
of the graph $\widehat{G}$ by adding to the graph
$\widehat{G}|_{\{\widehat{x}_1,...,\widehat{x}_{k-3}\}}$ its
$(G_1,H_1)$-extensions by the following rule. Consider all pairs
of $(G_1,H_1)$-extensions $(A,B)$ of
$\widehat{G}|_{\{\widehat{x}_1,...,\widehat{x}_{k-3}\}}$ in
$\widehat{G}$ such that $(E(A)\setminus
E(\widehat{G}|_{\{\widehat{x}_1,...,\widehat{x}_{k-3}\}}))\cap
(E(B)\setminus
E(\widehat{G}|_{\{\widehat{x}_1,...,\widehat{x}_{k-3}\}}))\neq\varnothing$.
If such pairs exist (we say that extensions from such pairs are
{\it intersecting}) we take their union and denote it by
$\widehat{X}^1_{\widehat{G}}(\widehat{x}_1,...,\widehat{x}_{k-3})$.

Suppose that the number of all the $(G_1,H_1)$-extensions such
that there exist other $(G_1,H_1)$-extensions which intersect them
is greater than $2(k-3)(k-2)$. Let us prove that
$\rho(\widehat{X}^1_{\widehat{G}}(\widehat{x}_1,...,\widehat{x}_{k-3}))>k-2$.

Let us reconsider the construction of the union
$\widehat{X}^1_{\widehat{G}}(\widehat{x}_1,...,\widehat{x}_{k-3})$.
Here we assume that the extensions are added step by step. At each
step we add either an extension intersecting extensions added
earlier or an intersecting pair of extensions which does not
intersect extensions added earlier. Let $v_i$ be a number of
vertices added at the $i$-th step. Let $h>(k-3)(k-2)$ be a number
of steps. Then
$$
 \frac{(k-2)(v_1+...+v_h)+h}{v_1+...+v_h+k-3}>k-2.
$$
Therefore,
$\rho(\widehat{X}_{\widehat{G}}^1(\widehat{x}_1,...,\widehat{x}_{k-3}))>k-2$.
Thus, the number of intersecting $(G_1,H_1)$-extensions is not
greater than $2(k-3)(k-2)$.

If in the graph $\widehat{G}$ there is no intersecting
$(G_1,H_1)$-extensions of
$\widehat{G}|_{\{\widehat{x}_1,...,\widehat{x}_{k-3}\}}$ we set
$$
 \widehat{X}_{\widehat{G}}^1(\widehat{x}_1,...,\widehat{x}_{k-3})=\widehat{G}|_{\{\widehat{x}_1,...,\widehat{x}_{k-3}\}}.
$$
Let
$$
 X_{\widehat{G}}^1(\widehat{x}_1,...,\widehat{x}_{k-3})=
 [\widehat{X}_{\widehat{G}}^1(\widehat{x}_1,...,\widehat{x}_{k-3})]_{\widehat{G}}.
$$

Consider a set $\mathcal{M}$ of all pairs
$([M]_{\widehat{G}},\widehat{G}|_{\{\widehat{x}_1,...,\widehat{x}_{k-3}\}}),$
where $M$ is an $(G_1,H_1)$-extension of the graph
$\widehat{G}|_{\{\widehat{x}_1,...,\widehat{x}_{k-3}\}}$ in
$\widehat{G}$ without vertices of the graph
$X_{\widehat{G}}^1(\widehat{x}_1,...,\widehat{x}_{k-3})\setminus
\widehat{G}|_{\{\widehat{x}_1,...,\widehat{x}_{k-3}\}}$. Choose
from the set $\mathcal{M}$ a collection of non-isomorphic pairs
$$
 ([M_1]_{\widehat{G}},\widehat{G}|_{\{\widehat{x}_1,...,\widehat{x}_{k-3}\}}),...,
 ([M_{\tau}]_{\widehat{G}},\widehat{G}|_{\{\widehat{x}_1,...,\widehat{x}_{k-3}\}}),
$$
such that $\tau$ is maximal (we call some pairs of graphs
$(A_1,B_1),(A_2,B_2)$, with $V(A_1)=\{a^1_1,...,a^1_n\}$,
$V(A_2)=\{a^2_1,...,a^2_n\}$, $V(B_1)=\{a^1_1,...,a^1_m\}$,
$V(B_2)=\{a^2_1,...,a^2_m\}$, $m<n$, {\it isomorphic} if $
 \{a^1_i,a^1_j\}\in E(A_1)\setminus E(B_1)\Leftrightarrow
 \{a^2_i,a^2_j\}\in E(A_2)\setminus E(B_2))$. Let
$$
 X_{\widehat{G}}^2(\widehat{x}_1,...,\widehat{x}_{k-3}) =
 [M_1]_{\widehat{G}}\cup...\cup[M_{\tau}]_{\widehat{G}}\cup
 X_{\widehat{G}}^1(\widehat{x}_1,...,\widehat{x}_{k-3}),
$$
$$
 \widehat{X}_{\widehat{G}}^2(\widehat{x}_1,...,\widehat{x}_{k-3})=
 M_1 \cup...\cup M_{\tau} \cup
 \widehat{X}_{\widehat{G}}^1(\widehat{x}_1,...,\widehat{x}_{k-3}).
$$

Similarly to the case of intersecting $(G_1,H_1)$-extensions, the
number of extenders of strict $(G_4^1,H_1)$- and
$(G_4^2,H_1)$-extensions of
$\widehat{G}|_{\{\widehat{x}_1,...,\widehat{x}_{k-3}\}}$ in
$\widehat{G}$ which have common vertices with other such extenders
is less than or equal to $2(k-3)(k-2)$. Let $W_1$ be the union of
all such extensions. Set $W_2=[W_1]_{\widehat{G}}$,
$$
 X_{\widehat{G}}^3(\widehat{x}_1,...,\widehat{x}_{k-3})=W_2\cup
 X_{\widehat{G}}^2(\widehat{x}_1,...,\widehat{x}_{k-3}),\,\,\,
 \widehat{X}_{\widehat{G}}^3(\widehat{x}_1,...,\widehat{x}_{k-3})=W_1\cup
 \widehat{X}_{\widehat{G}}^2(\widehat{x}_1,...,\widehat{x}_{k-3}).
$$

Consider a set $\mathcal{M}$ of all pairs
$([M]_{\widehat{G}},\widehat{G}|_{\{\widehat{x}_1,...,\widehat{x}_{k-3}\}}),$
where $M$ is $(G_4^1,H_1)$- or $(G_4^2,H_1)$-extension of the
graph $\widehat{G}|_{\{\widehat{x}_1,...,\widehat{x}_{k-3}\}}$ in
$\widehat{G}\setminus
(W_2\setminus\widehat{G}|_{\{\widehat{x}_1,...,\widehat{x}_{k-3}\}})$
having no vertices from
$V(X_{\widehat{G}}^3(\widehat{x}_1,...,\widehat{x}_{k-3})\setminus
\widehat{G}|_{\{\widehat{x}_1,...,\widehat{x}_{k-3}\}})$. Choose
 non-isomorphic pairs
$$
([M_1]_{\widehat{G}},\widehat{G}|_{\{\widehat{x}_1,...,\widehat{x}_{k-3}\}}),...,
([M_{\tau}]_{\widehat{G}},\widehat{G}|_{\{\widehat{x}_1,...,\widehat{x}_{k-3}\}})
$$
from the set $\mathcal{M}$ in such a way that the number $\tau$ is
maximal. Set
$$
 X_{\widehat{G}}^4(\widehat{x}_1,...,\widehat{x}_{k-3})=
 X_{\widehat{G}}^3(\widehat{x}_1,...,\widehat{x}_{k-3})
 \cup[M_1]_{\widehat{G}}\cup...\cup[M_{\tau}]_{\widehat{G}},
$$
$$
\widehat{X}_{\widehat{G}}^4(\widehat{x}_1,...,\widehat{x}_{k-3})=
\widehat{X}_{\widehat{G}}^3(\widehat{x}_1,...,\widehat{x}_{k-3})
\cup M_1 \cup...\cup M_{\tau}.
$$

For each vertex $\widehat{x}\in V(\widehat{G})\setminus
V(X_{\widehat{G}}^4(\widehat{x}_1,...,\widehat{x}_{k-3}))$
adjacent to less than or equal to $k-5$ vertices of
$\widehat{x}_1,...,\widehat{x}_{k-3}$ consider a set
$\Upsilon(\widehat{x})$ containing $\widehat{x}$ and all vertices
$\widehat{x}^1\in V(\widehat{G})\setminus
V(X_{\widehat{G}}^4(\widehat{x}_1,...,\widehat{x}_{k-3}))$,
satisfying the following property. There exists a vertex
$\widehat{x}^2\in V(\widehat{G})\setminus
V(X_{\widehat{G}}^4(\widehat{x}_1,...,\widehat{x}_{k-3}))$ such
that $\widehat{x}^1\sim\widehat{x}^2$,
$\widehat{x}^1\sim\widehat{x}$, $\widehat{x}^2\sim\widehat{x}$,
$\widehat{x}^1\sim\widehat{x}_i$,
$\widehat{x}^2\sim\widehat{x}_i$, $i\in\{1,...,k-3\}$. Let
$\Upsilon$ be a union of sets $\Upsilon(\widehat{x})$ over all
$\widehat{x}$ such that
$$
 f(X_{\widehat{G}}^4(\widehat{x}_1,...,\widehat{x}_{k-3})\cup\widehat{G}|_{\Upsilon(\widehat{x})},
 X_{\widehat{G}}^4(\widehat{x}_1,...,\widehat{x}_{k-3}))<0.
$$
Consider also a set $\mathcal{M}$ of all pairs
$(X_{\widehat{G}}^4(\widehat{x}_1,...,\widehat{x}_{k-3})\cup
[\widehat{G}|_{\Upsilon(\widehat{x})\cup\{\widehat{x}_1,...,\widehat{x}_{k-3}\}}]_{\widehat{G}},
X_{\widehat{G}}^4(\widehat{x}_1,...,\widehat{x}_{k-3}))$, where
$$
 f(X_{\widehat{G}}^4(\widehat{x}_1,...,\widehat{x}_{k-3})\cup\widehat{G}|_{\Upsilon(\widehat{x})},
 X_{\widehat{G}}^4(\widehat{x}_1,...,\widehat{x}_{k-3}))=0.
$$
Let
$([M_1]_{\widehat{G}},X_{\widehat{G}}^4(\widehat{x}_1,...,\widehat{x}_{k-3})),...,
([M_{\tau}]_{\widehat{G}},X_{\widehat{G}}^4(\widehat{x}_1,...,\widehat{x}_{k-3}))$
be non-isomorphic pairs from $\mathcal{M}$ such that the number
$\tau$ is maximal. Set
$$
 X_{\widehat{G}}^5(\widehat{x}_1,...,\widehat{x}_{k-3}) =
 [M_1]_{\widehat{G}}\cup...\cup[M_{\tau}]_{\widehat{G}}\cup
 X_{\widehat{G}}^4(\widehat{x}_1,...,\widehat{x}_{k-3})\cup
 [\widehat{G}|_{\Upsilon\cup\{\widehat{x}_1,...,\widehat{x}_{k-3}\}}]_{\widehat{G}},
$$
$$
 \widehat{X}_{\widehat{G}}^5(\widehat{x}_1,...,\widehat{x}_{k-3}) =
 M_1 \cup...\cup M_{\tau} \cup
 \widehat{X}_{\widehat{G}}^4(\widehat{x}_1,...,\widehat{x}_{k-3})\cup\widehat{G}|_{\Upsilon}
 \widehat{X}_{\widehat{G}}^5(\widehat{x}_1,...,\widehat{x}_{k-3}).
$$

\subsection{Graphs $X^l_j(\widehat{x}_1)$,
$\widehat{X}^l_j(\widehat{x}_1)$, $l\in\{1,2,3,4,5\}$}

Let $\widehat{G}_1,\widehat{G}_2,...$ be graphs satisfying the
following properties:

\begin{itemize}

\item[---]
$\bigcap_{i=1}^{\infty}V(\widehat{G}_i)\supseteq\{\widehat{x}_1,...,\widehat{x}_{k-3}\}$;

\item[---] pairs
$(\widehat{G}_1,\widehat{G}_1|_{\{\widehat{x}_1,...,\widehat{x}_{k-3}\}})$,
$(\widehat{G}_2,\widehat{G}_2|_{\{\widehat{x}_1,...,\widehat{x}_{k-3}\}})$,
... are pairwise non-isomorphic;

\item[---]
$\rho(X^5_{\widehat{G}_i}(\widehat{x}_1,...,\widehat{x}_{k-3}))\leq
k-2$ for every $i\in\mathbb{N}$;

\item[---] there is no such a graph $\widehat{G}_0$ that the
graphs $\widehat{G}_0,\widehat{G}_1,\widehat{G}_2,...$ satisfy the
first three properties.

\end{itemize}

Then obviously there exists a final set $\{i_1,...,i_{a(k)}\}$
such that pairs
$$
 (X^5_{\widehat{G}_{i_j}}(\widehat{x}_1,...,\widehat{x}_{k-3}),\widehat{G}_{i_j}|_{\{\widehat{x}_1,...,\widehat{x}_{k-3}\}})),
 \,\,\,
 j\in\{1,...,a(k)\}
$$ are pairwise non-isomorphic. Moreover,
for any $i\in\mathbb{N}$ there exists $j\in\{1,...,a(k)\}$ such
that the pairs
$(X^5_{\widehat{G}_{i}}(\widehat{x}_1,...,\widehat{x}_{k-3}),\widehat{G}_i|_{\{\widehat{x}_1,...,\widehat{x}_{k-3}\}}))$
and
$(X^5_{\widehat{G}_{i_j}}(\widehat{x}_1,...,\widehat{x}_{k-3}),\widehat{G}_{i_j}|_{\{\widehat{x}_1,...,\widehat{x}_{k-3}\}}))$
are isomorphic.

For any $j\in\{1,...,a(k)\}$ and $l\in\{1,2,3,4,5\}$ set
$$
 X^l_j(\widehat{x}_1,...,\widehat{x}_{k-3})=
 X^l_{\widehat{G}_{i_j}}(\widehat{x}_1,...,\widehat{x}_{k-3}).
$$

Consider the graph $\widehat{G}$. Let $\xi\in\{1,...,k-4\}$ be
fixed. Let $Y^1,...,Y^t$ be subgraphs of the subgraph
$\widehat{G}$ satisfying the following properties.
\begin{itemize}
\item[---] For each $i\in\{1,...,t\}$ there exist
$j\in\{1,...,a(k)\}$ and vertices
$\widehat{x}_{\xi+1}^i,...,\widehat{x}_{k-3}^i$ of the graph
$\widehat{G}$ such that
$Y^i=X^5_j(\widehat{x}_1,...,\widehat{x}_{\xi},\widehat{x}_{\xi+1}^i,...,\widehat{x}_{k-3}^i)$.

\item[---] Take arbitrary $i_1,i_2\in\{1,...,t\},$ $i_1\neq i_2,$
such that $\widehat{x}_i^{i_1},\widehat{x}_i^{i_2}$ do not
coincide for some $i\in\{\xi+1,...,k-4\}$. Let
$\mu\in\{\xi,...,k-4\}$ be such that for all
$i\in\{\xi+1,...,\mu\}$ the vertices
$\widehat{x}_i^{i_1},\widehat{x}_i^{i_2}$ coincide, and the
vertices $\widehat{x}^{i_1}_{\mu+1}$, $\widehat{x}^{i_2}_{\mu+1}$
are different. Then $f\left(\bigcup_{i\in
I^{i_1}}Y^i,\bigcap_{i\in I^{i_1}}Y^i\right)\leq 0$, where
$I^{i_1}=\{u:\,\forall i\in\{\xi+1,...,\mu\} \,\,
\widehat{x}_i^{u}=\widehat{x}_i^{i_1},\widehat{x}_{\mu+1}^{i_1}\neq\widehat{x}_{\mu+1}^{u}\}$.

\item[---]  Pairs $\left(\bigcup_{i\in I^{i_1}_1} Y^i,
\widehat{G}|_{\left\{\widehat{x}_1,...,\widehat{x}_{\xi},\widehat{x}^{i_1}_{\xi+1},...,\widehat{x}^{i_1}_{\mu}\right\}}\right)$,
$\left(\bigcup_{i\in I^{i_2}_2} Y^i,
\widehat{G}|_{\left\{\widehat{x}_1,...,\widehat{x}_{\xi},\widehat{x}^{i_1}_{\xi+1},...,\widehat{x}^{i_1}_{\mu}\right\}}\right)$
are non-isomorphic for any $\mu\in\{\xi,...,k-4\}$ and any
$i_1,i_2\in\{1,...,t\}$ such that
$\widehat{x}^{i_1}_{\mu+1}\neq\widehat{x}^{i_2}_{\mu+1}$. Here
$I_1^{i_1}=\{u: \, \forall i\in\{\xi+1,...,\mu+1\} \,\,
\widehat{x}^{i_1}_i=\widehat{x}^u_i\}$, $I_2^{i_2}=\{u: \, \forall
i\in\{\xi+1,...,\mu\} \,\,
\widehat{x}^{i_1}_i=\widehat{x}^j_i,\widehat{x}^{i_2}_{\mu+1}=\widehat{x}^u_{\mu+1}\}$.
\end{itemize}

If in the graph $\widehat{G}$ there is no subgraph $Y^{t+1}$
different from $Y^1,...,Y^t$ and such that $Y^1,...,Y^{t+1}$
satisfy the three properties described above, then we denote
$Y^1\cup...\cup Y^t$ by
$X^5_{\widehat{G}}(\widehat{x}_1,...,\widehat{x}_{\xi})$. If $t=0$
set
$X^5_{\widehat{G}}(\widehat{x}_1,...,\widehat{x}_{\xi})=\widehat{G}|_{\{\widehat{x}_1,...,\widehat{x}_{\xi}\}}$.

If, in addition, graphs $Y^1,...,Y^t$ follow the properties
described below then we say that the graph
$X^5_{\widehat{G}}(\widehat{x}_1,...,\widehat{x}_{\xi})$ is {\it
$(\widehat{x}_1,...,\widehat{x}_{\xi})$-net in $\widehat{G}$}.

\begin{itemize}

\item[---] For any $i_1,i_2\in\{1,...,t\},$ $i_1\neq i_2$, the set
$V(Y^1)\cap V(Y^2)$ is a subset of
$\{\widehat{x}_1,...,\widehat{x}_{\xi},\widehat{x}_{\xi+1}^{i_1},...,\widehat{x}_{k-4}^{i_1}\}\cap
\{\widehat{x}_1,...,\widehat{x}_{\xi},\widehat{x}_{\xi+1}^{i_2},...,\widehat{x}_{k-4}^{i_2}\}$.

\item[---] For any $i_1,i_2\in\{1,...,t\},$ $i_1\neq i_2$, either
the sets
$\{\widehat{x}^{i_1}_{\xi+1},...,\widehat{x}^{i_1}_{k-3}\},
\{\widehat{x}^{i_2}_{\xi+1},...,\widehat{x}^{i_2}_{k-3}\}$ do not
intersect each other, or there exists $\mu\in\{\xi+1,...,k-4\}$
such that for any $i\in\{\xi+1,...,\mu\}$ the vertices
$\widehat{x}_i^{i_1},\widehat{x}_i^{i_2}$ coincide and the sets
$\{\widehat{x}^{i_1}_{\mu+1},...,\widehat{x}^{i_1}_{k-3}\},
\{\widehat{x}^{i_2}_{\mu+1},...,\widehat{x}^{i_2}_{k-3}\}$ do not
intersect each other.

\item[---] For any $i_1,i_2\in\{1,...,t\},$ $i_1\neq i_2$, the set
$E(Y_{i_1}\cup Y_{i_2})\setminus(E(Y_{i_1})\cup E(Y_{i_2}))$ is
empty.

\end{itemize}

Let us consider the following situation. There exist graphs
$\widehat{G}^1,\widehat{G}^2$ and vertices
$\widehat{x}^1_1,...,\widehat{x}^1_{\xi},\widehat{x}^2_1,...,\widehat{x}^2_{\xi}$
such that

\begin{itemize}

\item[---]
$X^5_{\widehat{G}^1}(\widehat{x}^1_1,...,\widehat{x}^1_{\xi})=Y^1_1\cup...\cup
Y^t_1$,
$X^5_{\widehat{G}^2}(\widehat{x}^2_1,...,\widehat{x}^2_{\xi})=Y^1_2\cup...\cup
Y^t_2$;

\item[---] for any $i\in\{1,...,t\}$ there exists
$j\in\{1,...,a(k)\}$ such that
$Y^i_1=X^5_j(\widehat{x}^1_1,...,\widehat{x}^1_{\xi},\widehat{x}_{\xi+1}^{i,1},...,\widehat{x}_{k-3}^{i,1})$,
$Y^i_2=X^5_j(\widehat{x}^2_1,...,\widehat{x}^2_{\xi},\widehat{x}_{\xi+1}^{i,2},...,\widehat{x}_{k-3}^{i,2})$
for some vertices
$\widehat{x}_{\xi+1}^{i,1},...,\widehat{x}_{k-3}^{i,1},\widehat{x}_{\xi+1}^{i,2},...,\widehat{x}_{k-3}^{i,2}$;

\item[---] for any different $i_1,i_2\in\{1,...,t\}$ there exists
$\mu\in\{\xi,...,k-4\}$ such that
$$
 \widehat{x}_{j+1}^{i_1,1}=\widehat{x}_{j+1}^{i_2,1}, \,\,
 \widehat{x}_{j+1}^{i_1,2}=\widehat{x}_{j+1}^{i_2,2}, \,\,\,
 j\in\{\xi,...,\mu-1\},
$$
$$
 \widehat{x}_{\mu+1}^{i_1,1}\neq\widehat{x}_{\mu+1}^{i_2,1}, \,\,
 \widehat{x}_{\mu+1}^{i_1,2}\neq\widehat{x}_{\mu+1}^{i_2,2};
$$

\item[---] the graph
$X^5_{\widehat{G}^2}(\widehat{x}^2_1,...,\widehat{x}^2_{\xi})$ is
an $(\widehat{x}^2_1,...,\widehat{x}^2_{\xi})$-net in
$\widehat{G}^2$.

\end{itemize}

In this case we say that the graph
$X^5_{\widehat{G}^2}(\widehat{x}^2_1,...,\widehat{x}^2_{\xi})$ is
a {\it net of the graph
$X^5_{\widehat{G}^1}(\widehat{x}^1_1,...,\widehat{x}^1_{\xi})$}.
For any $i\in\{1,...,t\}$, $\mu\in\{\xi+1,...,k-3\}$ we introduce
a notation
$$
 \widehat{x}^{i,2}_{\mu}=\mbox{NET}_{\widehat{G}^1,\widehat{G}^2,
 \widehat{x}^1_1,...,\widehat{x}^1_{\xi},\widehat{x}^2_1,...,\widehat{x}^2_{\xi}}(\widehat{x}^{i,1}_{\mu}).
$$

Note that the function $\mbox{NET}$ is defined on the set of
symbols $\{\widehat{x}^{i,1}_{\mu},\,
i\in\{1,...,t\},\mu\in\{\xi+1,...,k-3\}\}$ of cardinality
$t(k-3-\xi)$. It means that some vertices from the set of vertices
$\{\widehat{x}^{i,1}_{\mu},\,
i\in\{1,...,t\},\mu\in\{\xi+1,...,k-3\}\}$ can be equal, but the
corresponding symbols are different. In this cases the function
$\mbox{NET}$ assigns a vertex to several different vertices.\\

Obviously the number of different $(\widehat{x}_1)$-nets (up to
isomorphism) in $\widehat{G}$ for different graphs $\widehat{G}$
with {\it maximal density} greater than or equal to $k-2$ is
finite (maximal density of a graph $G$ equals
$\rho^{\max}(G):=\max_{H\subseteq G}\{\rho(H)\}$). Let
$X^5_1(\widehat{x}_1),...,X^5_{\widetilde{m}(k)}(\widehat{x}_1)$
be such nets with the following property. For any net
$X^5_0(\widehat{x}_i)$ in $\widehat{G}$ with maximal density
greater than or equal to $k-2$ there exists a number
$j\in\{1,...,\widetilde{m}(k)\}$ such that the mapping
$X^5_0(\widehat{x}_i)\rightarrow X^5_j(\widehat{x}_i)$ preserving
the vertex order is an isomorphism. Let graphs
$X^5_1(\widehat{x}_1),...,X^5_{\widetilde{m}(k)}(\widehat{x}_1)$
be ordered so that maximal densities of the graphs
$X^5_1(\widehat{x}_1),...,X^5_{\widehat{m}(k)}(\widehat{x}_1)$ are
equal to $k-2$ and the densities of graphs
$X^5_1(\widehat{x}_1),...,X^5_{m(k)}(\widehat{x}_1)$ are equal to
$k-2$. For each number $j\in\{1,...,\widetilde{m}(k)\}$ let
$X^5_j(\widehat{x}_1)=Y^1_j\cup...\cup Y^{t(j)}_j$ be a
decomposition defined by properties described above (definition of
an $(\widehat{x}_1)$-net). Let $j\in\{1,...,\widetilde{m}(k)\}$,
$s\in\{1,...,t(j)\}$. Denote by $i^{j,s}$ a number from
$\{1,...,a(k)\}$ such that
$Y^{s}_j=X^5_{i^{j,s}}(\widehat{x}_1,\widehat{x}^{j,s}_2,...,\widehat{x}^{j,s}_{k-3})$
for some vertices
$\widehat{x}^{j,s}_2,...,\widehat{x}^{j,s}_{k-3}$. For every
$l\in\{1,2,3,4\}$ set
$$
 X_j^l(\widehat{x}_1)=X^l_{i^{j,1}}(\widehat{x}_1,\widehat{x}^{j,1}_2,...,\widehat{x}^{j,1}_{k-3})\cup...\cup
 X^l_{i^{j,t(j)}}(\widehat{x}_1,\widehat{x}^{j,t(j)}_2,...,\widehat{x}^{j,t(j)}_{k-3}),
$$
$$
 \widehat{X}_j^l(\widehat{x}_1)=\widehat{X}^l_{i^{j,1}}(\widehat{x}_1,\widehat{x}^{j,1}_2,...,\widehat{x}^{j,1}_{k-3})\cup...\cup
 \widehat{X}^l_{i^{j,t(j)}}(\widehat{x}_1,\widehat{x}^{j,t(j)}_2,...,\widehat{x}^{j,t(j)}_{k-3}),
$$
$$
 \widehat{X}_j^5(\widehat{x}_1)=\widehat{X}^5_{i^{j,1}}(\widehat{x}_1,\widehat{x}^{j,1}_2,...,\widehat{x}^{j,1}_{k-3})\cup...\cup
 \widehat{X}^5_{i^{j,t(j)}}(\widehat{x}_1,\widehat{x}^{j,t(j)}_2,...,\widehat{x}^{j,t(j)}_{k-3}).
$$

The definition of the graph
$X^5_{\widehat{G}}(\widehat{x}_1,...,\widehat{x}_{\xi})$ implies
that there exists an analogous decomposition:
$$
 X_{\widehat{G}}^5(\widehat{x}_1,...,\widehat{x}_{\xi})=
 \bigcup\limits_{i=1}^{t(\widehat{x}_1,...,\widehat{x}_{\xi},\widehat{G})}
 X^5_{\widehat{G}}(\widehat{x}_1,...,\widehat{x}_{\xi},\widehat{x}^{i}_{\xi+1}(\widehat{x}_1,...,\widehat{x}_{\xi}),
 ...,\widehat{x}^{i}_{k-3}(\widehat{x}_1,...,\widehat{x}_{\xi})).
$$

Note that any graph $G$ with maximal density $\rho$ has subgraphs
$H_1,H_2,$ $H_1\subseteq H_2$ such that
$\rho(H_1)=\rho(H_2)=\rho$, the graph $H_1$ is strictly balanced,
the pair $(H_2,H_1)$ is $1/\rho$-neutral chain, and either the
pair $(G,H_2)$ is $1/\rho$-safe or $H_2=G$ (a pair $(H_2,H_1)$ is
called {\it $\alpha$-neutral chain} if $H_2\supseteq H_1$ and
there exist graphs $K_1,...,K_r,T_1,...,T_{r-1}$ with the
following properties: $H_1=K_1\subset K_2\subset ... \subset
K_r=H_2$; $T_i\subset K_i$, $i\in\{1,...,r-1\}$; pairs
$((K_i\setminus K_{i-1})\cup T_{i-1}, T_{i-1})$, $i\in\{2,...,r\}$
are $\alpha$-neutral; for any $i\in\{2,...,r\}$ there are no edges
connecting vertices of the graph $K_i\setminus K_{i-1}$ and
vertices of the graph $K_{i-1}\setminus T_{i-1}$).

Denote by $X_j^*(\widehat{x}_1)$, $X_j^{**}(\widehat{x}_1)$ the
corresponding subgraphs of the graph $X_j^5(\widehat{x}_1)$ for
each $j\in\{1,...,\widetilde{m}(k)\}$. The graph
$X_j^*(\widehat{x}_1)$ is strictly balanced with density
$\rho^{\max}(X_j^5(\widehat{x}_1))$. The pair
$(X_j^{**}(\widehat{x}_1),X_j^{*}(\widehat{x}_1))$ is
$1/\rho^{\max}(X_j^5(\widehat{x}_1))$-neutral chain. Note that
when $j\in\{1,...,\widehat{m}(k)\}$ the graphs
$X_j^{**}(\widehat{x}_1)$, $X_j^5(\widehat{x}_1)$ are equal.

Let us prove that the graph $X_j^{*}(\widehat{x}_1)$ contains the
vertex $\widehat{x}_1$ if $j\in\{1,...,\widehat{m}(k)\}$. It is
easy to see that in the graph $X_j^*(\widehat{x}_1)$ there are at
least $k-1$ vertices such that any vertex $\widehat{x}^0$ of them
follows the property described below. There exists a number
$l\in\{1,...,t(j)\}$ such that the vertices
$\widehat{x}_1,\widehat{x}^{j,l}_2,...,\widehat{x}^{j,l}_{k-3}$
are adjacent to $\widehat{x}^0$ in the graph
$X_j^5(\widehat{x}_1)$ and $\widehat{x}^0\in
V(X^5_{i^{j,l}}(\widehat{x}_1,\widehat{x}^{j,l}_2,...,\widehat{x}^{j,l}_{k-3}))$.
Otherwise
$$
 \frac{C_{y_1}^2+...+C_{y_v}^2+v(k-4)+(k-2)+(k-2)(k-3)+C_{k-3}^2}{2k-5+v}\geq
 k-2,
$$
$$
 y_1^2+...+y_v^2-4v+(k-2)+2(k-2)(k-3)+(k-3)(k-4)\geq
 2(k-2)(2k-5)
$$
for some natural numbers $v,y_1,...,y_v$ such that
$y_1+...+y_v=k-2$.

A function $\phi(y_1,...,y_{k-2})=y_1^2+...+y_{k-2}^2-4|\{i:\,
y_i\neq 0\}|$ achieves its maximal value on the set
$\mathbb{Z}_+^{k-2}\cap\{y_1+...+y_{k-2}=k-2\}$ when
$(y_1,...,y_{k-2})=(k-2,0,...,0)$. Therefore,
$$
  (k-2)^2-4+(k-2)+2(k-2)(k-3)+(k-3)(k-4)\geq
 2(k-2)(2k-5).
$$
A contradiction is obtained. Thus, the vertex $\widehat{x}_1$ is
adjacent to $k-1$ or more vertices of the graph
$X_j^*(\widehat{x}_1)$. If the vertex $\widehat{x}_1$ is not a
vertex of this graph, then the density of the graph
$X_j^1(\widehat{x}_1)|_{V(X_j^*(\widehat{x}_1))\cup\{\widehat{x}_1\}}$
is greater than the density of the graph $X_j^*(\widehat{x}_1)$.

\section{Ehrenfeucht game}

The main tool in proofs of zero-one laws for the first order
properties of the random graphs is a result proved by
A.~Ehrenfeucht in 1960 (see \cite{Ehren}). In this section we
formulate its particular case for graphs. First of all let us
define the Ehrenfeucht game on two graphs $G,H$ with $i$ rounds
(see \cite{Janson}, \cite{Alon}, \cite{Veresh}--\cite{Shelah},
\cite{zhuk_dan}, \cite{Ehren}--\cite{Zhuk2}). Let
$V(G)=\{x_1,...,x_n\},$ $V(H)=\{y_1,...,y_m\}$. At the
$\nu\mbox{-}$th step ($1 \leq \nu \leq i$) Spoiler chooses a
vertex from any graph. He chooses either a vertex $x_{j_{\nu}}\in
V(G)$ or a vertex $y_{j'_{\nu}}\in V(H)$. At the same round
Duplicator chooses a vertex from the other graph. Let Spoiler
choose the vertex $x_{j_{\mu}}\in V(G),$ $j_{\mu}=j_{\nu}$
($\nu<\mu$), at the $\mu\mbox{-}$th round. Duplicator must choose
the vertex $y_{j'_{\nu}}\in V(H)$. If at this round Spoiler
chooses a vertex $x_{j_{\mu}}\in V(G),$
$j_{\mu}\notin\{j_1,...,j_{\mu-1}\},$ then Duplicator must choose
a vertex $y_{j'_{\mu}}\in V(H)$ such that
$j'_{\mu}\notin\{j'_1,...,j'_{\mu-1}\}.$ If Duplicator cannot find
such a vertex then Spoiler wins. After the final round vertices
$x_{j_1},...,x_{j_{i}}\in V(G)$, $y_{j'_1},...,y_{j'_{i}}\in V(H)$
are chosen. Some of these vertices probably coincide. Choose
pairwise different vertices: $x_{h_1},...,x_{h_l};$
$y_{h'_1},...,y_{h'_l},$ $l \leq i.$ Duplicator wins if and only
if the corresponding subgraphs are isomorphic:
$$
 G|_{\{x_{h_1},...,x_{h_l}\}}\cong
 H|_{\{y_{h'_1},...,y_{h'_l}\}}.
$$

\begin{theorem} [\cite{Ehren}]
Let $G,H$ be two graphs. Let $i\in\mathbb{N}$ be some natural
number. Duplicator has a winning strategy in the game $EHR(G,H,i)$
if and only if for any first-order property $L$ expressed by a
formula with quantifier depth at most $i$ either $G$ and $H$
satisfy $L$ or $G$ and $H$ do not satisfy $L$. \label{ehren}
\end{theorem}

In the two following sections we state lemmas which we use in the
proof of Theorem \ref{main} (see Section~\ref{theorem_proof}). We
prove lemmas in Section~\ref{lemma_proof}.

\section{Main lemmas}

Let $\mathcal{G}\in\Omega_N,$ $\widetilde{x}\in\mathcal{V}_N$,
$\mathcal{G}\supset Y \supset \mathcal{G}|_{\{\widetilde{x}\}}$.
We call the pair $(Y,\widetilde{x})$ {\it $j$-maximal in
$\mathcal{G}$}, where $j\in\{1,...,\widehat{m}(k)\}$, if
$Y=X_j^5(\widetilde{x})=X^5_{\mathcal{G}}(\widetilde{x})$.


Let $\mathcal{L}^k_j(N)\subset\Omega_N$ be a set of graphs
$\mathcal{G}$ such that there are a vertex $\widetilde{x}$ and a
graph $Y$ such that the pair $(Y,\widetilde{x})$ is
$j$-maximal in $\mathcal{G}$. 
Set
$$
 \mathcal{L}^k_{\widehat{m}(k)+1}(N)=\Omega_N\setminus(\mathcal{L}^k_1(N)\cup...\cup\mathcal{L}^k_{\widehat{m}(k)}(N)),
$$
$$
 \mathcal{A}^k_{j_1...j_t}(N)=\left(\bigcap_{i=1}^t\mathcal{L}^k_{j_i}\right)\bigcap
 \left(\Omega_N\setminus\left(\bigcup_{i\in\{1,...,\widehat{m}(k)\}\setminus\{j_1,...,j_t\}}
 \mathcal{L}^k_{i}(N)\right)\right)
$$
for any different $j_1,...,j_t\in\{1,...,\widehat{m}(k)+1\}$. The
following lemma provides pairs of graphs
$(\widetilde{G},\widetilde{H})$ such that Duplicator has a winning
strategy in the game $EHR(\widetilde{G},\widetilde{H},k)$.

\begin{lemma}
For any subset $\{j_1,...,j_t\}\subset\{1,...,\widehat{m}(k)+1\}$
Duplicator has a winning strategy in the game
$EHR(\widetilde{G},\widetilde{H},k)$ for almost all pairs of
graphs $(\widetilde{G},\widetilde{H})$ from
$\mathcal{A}^k_{j_1...j_t}(N)\times\mathcal{A}^k_{j_1...j_t}(M)$.
\label{classes}
\end{lemma}

In the following lemma an asymptotic behavior of probabilities of
$\mathcal{A}^k_{j_1...j_t}(N)$ is described.

\begin{lemma}
For any different $j_1,...,j_t\in\{1,...,\widehat{m}(k)+1\}$ there
exist constants $0\leq\xi_{j_1...j_t} \leq 1$ such that
\begin{equation}
 \lim\limits_{N\rightarrow\infty}{\sf
 P}_{N,p}(\mathcal{A}^k_{j_1...j_t}(N))=
 \xi_{j_1...j_t}.
\label{A}
\end{equation}

\label{L}
\end{lemma}

\section{Auxiliary lemmas}
\label{additional_lemma}

In this section we give two statements which we use in the proofs
of Lemma \ref{classes} and Lemma \ref{L}.

\begin{lemma} \label{balance-uho}
Let $j_1,...,j_l\in\{1,...,\widehat{m}(k)\}$ be some numbers (some
of them may be equal). Let for any $i_1\in\{1,...,l\}$ there exist
$i_2\in\{1,...,l\}\setminus\{i_1\}$ such that the graphs
$X^{**}_{j_{i_1}}(\widetilde{x}^{i_1}),X^{**}_{j_{i_2}}(\widetilde{x}^{i_2})$
have a common vertex. Then
$\rho(X^{**}_{j_1}(\widetilde{x}^1)\cup...\cup
X^{**}_{j_l}(\widetilde{x}^l))=k-2$ if and only if the sets
$V(X^*_{j_i}(\widetilde{x}^i))$ coincide for all $i\in
\{1,...,l\}$. If not all of the sets coincide then
$\rho(X^{**}_{j_1}(\widetilde{x}^1)\cup...\cup
X^{**}_{j_l}(\widetilde{x}^l))>k-2$.
\end{lemma}

For an arbitrary graph $G$ we set $f(G)=v(G)-\alpha\cdot e(G)$.

\begin{lemma}
Let $G$ be a strictly balanced graph, $\rho(G)<k-2$. Let
$\mathcal{R}$ be the number of all $(K,T)$-maximal copies of the
graph $G$ in $G(N,p)$ for all $\alpha$-neutral pairs $(K,T)$ such
that $v(T)\leq k^3$, $v(K,T)\leq k^3$. Then $\mathcal{R}$
converges in probability to infinity. The fraction
$\frac{\mathcal{R}}{{\sf E}_{N,p}\mathcal{R}}$ converges in
probability to 1 and ${\sf E}_{N,p}\mathcal{R}=\Theta(N^{f(G)})$.
\label{max_graph}
\end{lemma}

We prove Lemma \ref{balance-uho} and Lemma \ref{max_graph} in
Section~\ref{lemma_proof}.

\section{Proof of Theorem \ref{main}}
\label{theorem_proof}

It follows from Lemma \ref{classes} that there exists a set
$\widetilde{\Omega}_N\subset\Omega_N,$ ${\sf
P}_{N,p}(\widetilde{\Omega}_N)\rightarrow 1,\,N\rightarrow\infty$,
and a partition of this set $\Omega_N^1,...,\Omega_N^{s(k)}$,
$\bigsqcup_{i=1}^{s(k)}\Omega_N^i=\widetilde{\Omega}_N$, such that
for any $i\in\{1,...,s(k)\},$ $N,M\in\mathbb{N}$ and any pair of
graphs $\mathcal{G}\in\widetilde{\Omega}_N^i$,
$\mathcal{H}\in\widetilde{\Omega}_M^i$ Duplicator has a winning
strategy in the game $EHR(\mathcal{G},\mathcal{H},k)$. Any set
$\Omega_N^i$ is an intersection of sets
$\mathcal{A}^k_{j_1...j_t}(N)$ for some
$j_1,...,j_t\in\{1,...,\widehat{m}(k)+1\}$ with
$\widetilde{\Omega}_N\subset\Omega_N$. Let $L$ be a first order
property expressed by a formula with a quantifier depth at most
$k$. By Theorem \ref{ehren} for each $i\in\{1,...,s(k)\}$ its
truth is the same for all graphs from
$\bigcup_{N\in\mathbb{N}}\Omega_N^i$. Lemma \ref{L} provides a
convergence of a probability of $\Omega_N^i$ for any
$i\in\{1,...,s(k)\}$. The subset
$\mathcal{A}_N(L)\subset\widetilde{\Omega}_N$ consisting of all
graphs satisfying the property $L$ is the union of $\Omega_N^i$,
$i\in I$, for some $I\subseteq \{1,...,s(k)\}$. Therefore, ${\sf
P}_{N,p}(\mathcal{A}_N(L))$ converges too. Theorem is proved.

\section{Proofs of lemmas}
\label{lemma_proof}

We do not give a proof of Lemma \ref{max_graph} in the paper
because it is a simple version of the proof of Theorem
\ref{maximal_extensions_neutral} that was proved in
\cite{zhuk_extensions}. The proof of Lemma \ref{L} is based on
Lemma \ref{balance-uho}. The proof of Lemma \ref{classes} is based
on Lemma \ref{max_graph}. Therefore, we prove Lemma
\ref{balance-uho} first and prove Lemma \ref{L}  and Lemma
\ref{classes} after that.

\subsection{Proof of Lemma
\ref{balance-uho}}

Consider some $\alpha$-neutral chain $(G,H)$ and graphs
$K_1,...,K_r,T_1,...,T_{r-1}$ such that

\begin{itemize}

\item[---] $H=K_1\subset K_2\subset ... \subset K_r=G$,
$T_i\subset K_i$, $i\in\{1,...,r-1\}$,

\item[---] the pairs $((K_i\setminus K_{i-1})\cup T_{i-1},
T_{i-1})$, $i\in\{2,...,r\}$, are $\alpha$-neutral,

\item[---] for any $i\in\{2,...,r\}$ the vertices of the graph
$K_i\setminus K_{i-1}$ and the vertices of the graph
$K_{i-1}\setminus T_{i-1}$ are not adjacent.

\end{itemize}

Suppose that $H$ is a strictly balanced graph with the density
$\rho(H)=1/\alpha$. Let us prove that the graph $G$ is balanced.\\

Let $F$ be a proper subgraph of $G$, $F_1=F\cap H,$
$F_i=F\cap(K_i\setminus K_{i-1}),$ $i\in\{2,...,r\}$. From the
definition of an $\alpha$-neutral pair it follows that
$$
 f(F_i\cup T_{i-1},T_{i-1})\geq 0,\,\,i\in\{2,...,r\}.
$$
Obviously,
$$
e(F_i\cup...\cup F_1,F_{i-1}\cup...\cup F_1)\leq e(F_i\cup
T_{i-1},T_{i-1}), \,\, i\in\{2,...,r\},
$$
$$
v(F_i\cup...\cup F_1,F_{i-1}\cup...\cup F_1)= v(F_i\cup
T_{i-1},T_{i-1}),\,\, i\in\{2,...,r\}.
$$
Therefore,
$$
f(F_i\cup...\cup F_1,F_{i-1}\cup...\cup F_1)\geq
0,\,\,i\in\{2,...,r\}. \label{balance_1}
$$
Furthermore,
$$
 f(F_1)\geq f(H)=0
$$
as $H$ is a strictly balanced graph. The last inequality is strict
if and only if $F_1\neq\varnothing$. Finally, we get
$$
 \rho(F)=(k-2)-\frac{(k-2)f(F)}{v(F)}=
$$
$$
 =(k-2)-\frac{(k-2)(f(F_1)+\sum_{i=2}^rf(F_i\cup...\cup F_1,
 F_{i-1}\cup...\cup F_1))}{v(F)}
 \leq k-2=\rho(G).
$$
Therefore, the graph $G$ is balanced.\\

For each $i\in\{1,...,l\}$ denote by $Y_i$ the graph
$X^{**}_{j_i}(\widetilde{x}_i)$. We prove Lemma \ref{balance-uho}
by induction. Consider the case $l=2$. Set $Y_{1,2}=Y_1\cap Y_2$.
Consider the following three situations.
\begin{itemize}
\item[1)] The set $V(X_{j_1}^*(\widetilde{x}^1))\cap V(Y_{1,2})$
is not empty. The graph $X_{j_1}^*(\widetilde{x}^1)\cap Y_{1,2}$
is a proper subgraph of the graph $X_{j_1}^*(\widetilde{x}^1)$.

\item[2)] The set $V(X_{j_1}^*(\widetilde{x}^1))\cap V(Y_{1,2})$
is empty.

\item[3)] The equality $X_{j_1}^*(\widetilde{x}^1)\cap
Y_{1,2}=X_{j_1}^*(\widetilde{x}^1)$ holds.
\end{itemize}

The graph $X_{j_1}^*(\widetilde{x}^1)$ is strictly balanced with
the density equal to $k-2$. The pair
$(Y_1,X_{j_1}^*(\widetilde{x}^1))$ ia an $1/(k-2)$-neutral chain.
Therefore, the graph $Y_1$ is balanced. Thus, in the first case
$$
 f(Y_1,X_{j_1}^*(\widetilde{x}^1)\cup Y_{1,2})\leq 0, \,\,\,\,\,
 f(Y_2)=0.
$$
Furthermore, the graph $X_{j_1}^*(\widetilde{x}^1)$ is strictly
balanced, $f(X_{j_1}^*(\widetilde{x}^1))=0$. Therefore,
$$
 \frac{e(X_{j_1}^*(\widetilde{x}^1))-e(X_{j_1}^*(\widetilde{x}^1),Y_{1,2}\cap X_{j_1}^*(\widetilde{x}^1))}
 {v(X_{j_1}^*(\widetilde{x}^1))-v(X_{j_1}^*(\widetilde{x}^1),Y_{1,2}\cap
 X_{j_1}^*(\widetilde{x}^1))}=\frac{e(Y_{1,2}\cap X_{j_1}^*(\widetilde{x}^1))}{v(Y_{1,2}\cap
 X_{j_1}^*(\widetilde{x}^1))}<k-2.
$$
So,
\begin{equation}
 f(X_{j_1}^*(\widetilde{x}^1)\cup Y_{1,2},Y_{1,2})\leq
 f(X_{j_1}^*(\widetilde{x}^1),Y_{1,2}\cap X_{j_1}^*(\widetilde{x}^1))< 0.
\label{strict}
\end{equation}
Finally, we get
$$
 e(Y_1\cup Y_2)\geq e(Y_2)+e(Y_1,X_{j_1}^*\cup Y_{1,2})+e(X_{j_1}^*(\widetilde{x}^1)\cup
 Y_{1,2},Y_{1,2})>(k-2)v(Y_1\cup
 Y_2).
$$
The last inequality is strict due to (\ref{strict}). Thus,
$$
 \rho(Y_1\cup
 Y_2)>k-2.
$$

Consider the second case: ($X_{j_1}^*(\widetilde{x}^1)\cap
Y_{1,2}=\varnothing$). From the definition of an $\alpha$-neutral
chain it follows that $v(Y_{1,2})-\alpha \cdot e(Y_{1,2})>0.$
Therefore,
$$
 e(Y_1\cup Y_2)\geq
 e(Y_2)+e(Y_1,Y_{1,2})=e(Y_{2})+e(Y_1)-e(Y_{1,2})>
$$
$$
 >
 (k-2)(v(Y_2)+v(Y_1)-v(Y_{1,2}))
 =(k-2)v(Y_1\cup Y_2).
$$
We get $\rho(Y_1\cup
Y_2)>k-2$.\\

Let finally $Y_{1,2}\supseteq X_{j_1}^*(\widetilde{x}^1).$ Then
$X_{j_1}^*(\widetilde{x}^1)=X_{j_2}^*(\widetilde{x}^2)$. Actually
if $X_{j_1}^*(\widetilde{x}^1)\cap
X_{j_2}^*(\widetilde{x}^2)\notin\{X_{j_1}^*(\widetilde{x}^1),X_{j_2}^*(\widetilde{x}^2),\varnothing\}$
then the pair $(X_{j_1}^*(\widetilde{x}^1)\cup
X_{j_2}^*(\widetilde{x}^2),X_{j_1}^*(\widetilde{x}^1))$ is
$\alpha$-rigid as the graph $X_{j_2}^*(\widetilde{x}^2)$ is
strictly balanced. This fact is in conflict with the properties
$X_{j_1}^*(\widetilde{x}^1)\cup
X_{j_2}^*(\widetilde{x}^2)\subseteq Y_2$ and
$\rho(Y_2)=\rho(X_{j_1}^*(\widetilde{x}^1))$ as the graph $Y_2$ is
balanced. In the cases $X_{j_1}^*(\widetilde{x}^1)\subset
X_{j_2}^*(\widetilde{x}^2)$, $X_{j_1}^*(\widetilde{x}^1)\subset
X_{j_2}^*(\widetilde{x}^2)$ we also get rigid pairs and obtain
contradiction. If $X_{j_1}^*(\widetilde{x}^1)\cap
X_{j_2}^*(\widetilde{x}^2)=\varnothing$ then the graph
$Y_2\setminus X_{j_2}^*(\widetilde{x}^2)$ contains the subgraph
$X_{j_1}^*(\widetilde{x}^1)$ with the density $k-2$. It is
impossible since
$(Y_2,X_{j_2}^*(\widetilde{x}^2))$ is $\alpha$-neutral chain.\\

Consider $l\geq 3$ pairs $(Y_i,\widetilde{x}^i)$. Let $V(Y_1)\cap
V(Y_2)\neq\varnothing$, $V(X_{j_1}^*(\widetilde{x}^1))\neq
V(X_{j_2}^*(\widetilde{x}^2))$,
$$
 \rho(Y_1\cup...\cup Y_{l-1})>k-2.
$$
The graph $\bigcup\limits_{i=1}^{l-1} Y_i\cap Y_l$ is a subgraph
of the graph $Y_l$. We have proved that $Y_l$ is a balanced graph.
Therefore,
$$
 \rho\left(\bigcup\limits_{i=1}^{l-1}Y_i\cap
 Y_l\right)\leq k-2.
$$
Thus,
$$
 \rho\left(\bigcup\limits_{i=1}^{l}Y_i\right)=
 \frac{e\left(\bigcup\limits_{i=1}^{l}Y_i\right)}
 {v\left(\bigcup\limits_{i=1}^{l}Y_i\right)}\geq
 \frac{e\left(\bigcup\limits_{i=1}^{l-1}Y_i\right)+e(Y_l)
 -e\left(\bigcup\limits_{i=1}^{l-1}Y_i\cap Y_l\right)}
 {v\left(\bigcup\limits_{i=1}^{l-1}Y_i\right)+v(Y_l)
 -v\left(\bigcup\limits_{i=1}^{l-1}Y_i\cap Y_l\right)}>k-2.
$$
So, the density equals $k-2$ if and only if for any graphs
$Y_{i_1},Y_{i_2}$ with common vertices
$X_{j_{i_1}}^*(\widetilde{x}^{i_1})=X_{j_{i_2}}^*(\widetilde{x}^{i_2})$.
For any $i_1\in\{1,...,l\}$ there exists
$i_2\in\{1,...,l\}\setminus\{i_1\}$ such that the graphs
$Y_{i_1},Y_{i_2}$ have a common vertex. Therefore, the density
equals $k-2$ if and only if the sets
$V(X^*_{j_i}(\widetilde{x}^i))$
coincide for all $i\in \{1,...,l\}$. Lemma is proved.\\

\subsection{Proof of Lemma   \ref{L}}

Let us prove the convergence of ${\sf
P}_{N,p}(\mathcal{L}_j^k(N))$ to some number $\xi_j$ for each
$j\in\{1,...,m(k)\}$ as $N\rightarrow\infty$. The proof is based
on three statements. The first one, Lemma \ref{balance-uho}, is
already proved. The second one is stated and proved in
\cite{Zhuk2}. An analogue of the third
statement is proved there too. Let us introduce some notation.\\

Let $j\in\{1,...,m(k)\}$. Let $v_j$ and $e_j$ be the numbers of
vertices and edges in the graph $X_j^5(\widehat{x}_1)$
respectively. Let $a_j$ be the number of automorphisms of the
graph $X_j^5(\widehat{x}_1)$ with the fixed point $\widehat{x}_1$.
Consider all ordered collections of $v_j$ vertices of the set
$\mathcal{V}_N$. Let us define a subset $M_j$ of the set of all
such collections.
\begin{itemize}

\item $M_j$ contains all different unordered collections.

\item Let $\left(\widetilde{x}_{i_1},...,
\widetilde{x}_{i_{v_j}}\right)\in M_j$. Let $\widetilde{Y}$ be a
graph on the set of vertices $\left\{\widetilde{x}_{i_1},...,
\widetilde{x}_{i_{v_j}}\right\}$. Assume that $\widetilde{Y}$ is a
strict $\left(X^5_j(\widehat{x}_1),
\widehat{G}|_{\{\widehat{x}_1\}}\right)$-extension of the graph
$\widetilde{Y}|_{\{\widetilde{x}_{i_1}\}}$. Let a graph obtained
by permutation of vertices
$\left(%
\begin{array}{ccc}
  i_{2} & ... & i_{v_j} \\
  t_{2} & ... & t_{v_j} \\
\end{array}%
\right)$ of the graph $\widetilde{Y}$ be a strict
$\left(X_j^5(\widehat{x}_1),\widehat{G}|_{\{\widehat{x}_1\}}\right)$-extension
of the graph $\widetilde{Y}|_{\{\widetilde{x}_{i_1}\}}$. Then the
set $M_j$ does not contain the collection
$\left(\widetilde{x}_{i_1},\widetilde{x}_{t_2},...,\widetilde{x}_{t_v}\right)$.
Otherwise this collection is in $M_j$.

\item In $M_j$ there are no collections except the described ones.

\end{itemize}

Set $m_j=|M_j|$. Let us enumerate all the collections from the set
$M_j$ by numbers $1,...,m_j$. Consider events
$\mathcal{B}_1^j,...,\mathcal{B}_{m_j}^j.$ The event
$\mathcal{B}_i^j$ is that a subgraph $\widetilde{Y}_i$ on the
$i$-th collection from $M_j$ and its first vertex form a
$j$-maximal pair. Let $A_i^j$ be an indicator of the event
$\mathcal{B}_i^j$. Consider a random variable
$A_j=\sum\limits_{i=1}^{m_j} A_i^j$ equal to a number of all
$j$-maximal pairs. We get
$$
 {\sf P}_{N,p}(A_j=0)=1-\sum_{i=1}^{m_j} {\sf P}_{N,p}(\mathcal{B}^j_i)+
 \sum_{i_1,i_2=1}^{m_j}{\sf P}_{N,p}(\mathcal{B}^j_{i_1}\cap\mathcal{B}^j_{i_2})+...+
$$
\begin{equation}
 +(-1)^n\sum_{i_1,i_2,...,i_n=1}^{m_j}{\sf P}_{N,p}
 (\mathcal{B}^j_{i_1}\cap\mathcal{B}^j_{i_2}\cap...\cap\mathcal{B}^j_{i_n})+...
\label{bruno1}
\end{equation}
The summation is over all different collections with  pairwise
different numbers. Let us prove that there exists a number $\xi_j$
such that
$$
 \lim\limits_{N\rightarrow\infty}{\sf P}_{N,p}(\mathcal{L}^k_j(N))=
 \lim\limits_{N\rightarrow\infty}(1-{\sf P}_{N,p}(A_j=0))=\xi_j.
$$

Let $\phi^j_1(N)$ be the probability that the pair
$(\widetilde{Y},\widetilde{x}_{i_1})$,
$\widetilde{Y}=\mathcal{G}|_{\{\widetilde{x}_{i_1},...,\widetilde{x}_{i_{v_j}}\}}$,
$\mathcal{G}\in\Omega_N$, is $j$-maximal under the condition that
the graph $\widetilde{Y}$ is a strict
$(X_j^5(\widehat{x}_1),\widehat{G}|_{\{\widehat{x}_1\}})$-extension
of the graph $\widetilde{Y}|_{\{\widetilde{x}_1\}}$. Then
$$
 \sum_{i=1}^{m_j} {\sf P}_{N,p}(\mathcal{B}^j_i)= {\sf E}_{N,p}(A_j) =
 NC_{N-1}^{v_j-1} \frac{(v_j-1)!}{a_j} \phi^j_1(N) p^{e_j} \sim
 \frac{\phi^j_1(N)}{a_j},
$$
where ${\sf E}_{N,p}(X_j)={\sf E}(X_j(G(N,p)))$. Set
$$
 a^j_n(N)=\sum_{i_1,i_2,...,i_n=1}^{m_j}{\sf P}_{N,p}
 (\mathcal{B}^j_{i_1}\cap\mathcal{B}^j_{i_2}\cap...\cap\mathcal{B}^j_{i_n}).
$$
We use the notation $i_1\sim i_2$ in the following case: the
numbers $i_1,i_2$ are from $\{1,...,m_j\}$, $i_1\neq i_2$, and the
collections from $M_j$ numerated by $i_1,i_2$ have common
vertices. Denote the sum with intersecting collections of vertices
by $r_j(n,N)$. In other words
$$
 a^j_n(N)
 -\sum_{i_1,i_2,...,i_n:\,\,\forall t_1\neq t_2\in\{1,...,n\} \,\, i_{t_1}\neq i_{t_2},i_{t_1}\nsim i_{t_2}}
 {\sf P}_{N,p}
 (\mathcal{B}^j_{i_1}\cap\mathcal{B}^j_{i_2}\cap...\cap\mathcal{B}^j_{i_n})=r_j(n,N).
$$

Let $\mathbf{Y}_1,...,\mathbf{Y}_n$ be pairwise disjoint
collections from $\mathcal{V}_N$ with cardinality $v_j$. Let
$\widetilde{x}^i_t$ be a vertex numerated by $t$ in the $i$-th
collection, $i\in\{1,...,n\}$, $t\in\{1,...,v_j\}$. Let
$\mathcal{Y}^j_n(N)$ be a set of all graphs $\mathcal{G}$ from
$\Omega_N$ such that for any $i\in\{1,...,n\}$ the pair
$\left(\mathcal{G}|_{\mathbf{Y}_i},\widetilde{x}^i_1\right)$ is
$j$-maximal in $\mathcal{G}$. Denote by $\mathcal{X}^j_n(N)$ a set
of all graphs $\mathcal{G}$ in $\Omega_N$ such that the subgraphs
$\mathcal{G}|_{\mathbf{Y}_1},...,\mathcal{G}|_{\mathbf{Y}_n}$ are
strict $\left(X^5_j(\widehat{x}_1),
\widehat{G}|_{\{\widehat{x}_1\}}\right)$-extensions of graphs
$\mathcal{G}|_{\{\widetilde{x}^1_1\}},...,
\mathcal{G}|_{\{\widetilde{x}^n_1\}}$ respectively. Let
$$
 \phi^j_n(N)={\sf
 P}_{N,p}(\mathcal{Y}^j_n(N)|\mathcal{X}^j_n(N)).
$$
Obviously the probability $\phi^j_n(N)$ does not depend on a
choice of sets $\mathbf{Y}_1,...,\mathbf{Y}_n$. We get
$$
 a^j_n(N)=\sum_{i_1,i_2,...,i_n=1}^{m_j}{\sf P}_{N,p}
 (\mathcal{B}^j_{i_1}\cap\mathcal{B}^j_{i_2}\cap...\cap\mathcal{B}^j_{i_n})
\sim \frac{\phi^j_n(N)}{n!}\left(\frac{1}{a_j}\right)^n+r_j(n,N).
$$

Let us formulate a statement from \cite{Zhuk2} (see Statement 3).

\begin{state}
Let $\{a_n(N)\}_{n\in\mathbb{N}}$ be a set of functions such that
there exists a sequence of numbers $\{b_n\}_{n\in\mathbb{N}}$
obeying the following law: $\forall n\in\mathbb{N}$
$a_n(N)\rightarrow b_n,$ $N\rightarrow\infty$. Let
$\sum\limits_{n=1}^{\infty}b_n=b$. If for any $N\in\mathbb{N}$ the
series $\sum\limits_{n=1}^{\infty}a_n(N)$ converges and for every
$s\in\mathbb{N}$, $N\in\mathbb{N}$
$$
 \sum\limits_{n=1}^{2s-1}a_n(N)\leq
 \sum\limits_{n=1}^{\infty}a_n(N)\leq
 \sum\limits_{n=1}^{2s}a_n(N),
$$
then $\sum\limits_{n=1}^{\infty}a_n(N)\rightarrow b$,
$N\rightarrow\infty$. \label{sum}
\end{state}

For any $n\in\mathbb{N}$, $N\in\mathbb{N}$ the inequality
$a^j_n(N)\geq a^j_{n+1}(N)$ holds. Therefore by Statement
\ref{sum} the convergence
$$
 {\sf P}_{N,p}(\mathcal{L}_j^k(N))\rightarrow \xi_j
$$
follows from the following fact. For each $n\in\mathbb{N}$
\begin{equation}
 \lim\limits_{N\rightarrow\infty}(a^j_n(N)-r_j(n,N))= b_j(n),
\label{lim_existense}
\end{equation}
\begin{equation}
 \lim\limits_{N\rightarrow\infty}r_j(n,N)=r_j(n),
\label{trash}
\end{equation}
\begin{equation}
 \sum_{n=1}^{\infty}(-1)^n(b_j(n)+r_j(n))<\infty.
\label{row}
\end{equation}

The equality (\ref{lim_existense}) follows from a statement
similar to Statement 2 from \cite{Zhuk2}. We do not give here a
proof of the statement because the proofs of these two statements
are the same.

\begin{state}
There exists $0<\zeta_j<1$ such that
$$
 \phi_1^j(N)\sim \zeta_j,
 \,\,\,
 \phi^j_n(N)\sim \zeta_j^n.
$$
\label{probab_converge}
\end{state}

All that remains is to show that equalities (\ref{trash}) and
$$
 \lim\limits_{n\rightarrow\infty} r_j(n)=0
$$
hold.

Let $x$ be a vertex. Let us define sets
$\mathcal{Q}^n_1(x)$,$\mathcal{Q}^n_2(x)$ in the following way:
$$
 (Q,x)\in\mathcal{Q}^n_1\Leftrightarrow
 ((v_j<v(Q)<nv_j)\wedge(\rho(Q)=k-2)\wedge(\exists
 \widehat{Y}_1,...,\widehat{Y}_n\,
 (\forall i\in\{1,...,n\}\
$$
$$
 ((\widehat{Y}_i,x)\cong
 (\widehat{Y}_j(\widehat{x}),\widehat{x}))\wedge
 (\forall i_1,i_2\in\{1,...,n\}\,\,(\widehat{Y}_{i_1}\cap\widehat{Y}_{i_2}\supset
 X_j^*(x)))
 \wedge(Q=\widehat{Y}_1\cup
 ...\cup \widehat{Y}_n))),
$$
$$
 (Q,x)\in\mathcal{Q}^n_2\Leftrightarrow
 ((v_j<v(Q)<nv_j)\wedge(\rho(Q)>k-2)\wedge(\exists
 \widehat{Y}_1,...,\widehat{Y}_n
$$
$$
 \forall i_1\in\{1,...,n\}\,\,
 ((\widehat{Y}_{i_1}\cong \widehat{Y}_j(\widehat{x}))\wedge(\exists i_2 \in\{1,...,n\} \,
 (\widehat{Y}_{i_1}\cap \widehat{Y}_{i_2}\neq\varnothing)))\wedge(Q=\widehat{Y}_1\cup
 ...\cup \widehat{Y}_n))).
$$
Let $(Q_1,x)\in\mathcal{Q}^{r_1}_{i_1}$, ...,
$(Q_t,x)\in\mathcal{Q}^{r_t}_{i_t}$, $i_1,...,i_t\in\{1,2\}$,
$r_1,...,r_t\in\mathbb{N}$, $r_1+...+r_t=n$. Let
$\widehat{Y}^l_1\cup...\cup \widehat{Y}^l_{r_l}$ be a
decomposition of $Q_l$, $l\in\{1,...,t\}$, into graphs isomorphic
to $X^5_j(\widehat{x})$. Let us introduce different collections of
vertices from $\mathcal{V}_N$ for graphs $Q_l$, $l\in\{1,...,t\}$,
in the same way as  $M_j$ was introduced. The first vertex in a
collection is fixed if and only if
$(Q_l,x)\in\mathcal{Q}^{r_l}_1$. For every $l\in\{1,...,t\}$
define an event $\mathcal{B}_i(Q_l)$. Its definition depends on
whether $(Q_l,x)\in\mathcal{Q}_1^{r_l}$ or
$(Q_l,x)\in\mathcal{Q}^{r_l}_2$. If
$(Q_l,x)\in\mathcal{Q}^{r_l}_2$ then the event
$\mathcal{B}_i(Q_l)$ is that the subgraph induced on the $i$-th
collection is isomorphic to $Q_l$.  If
$(Q_l,x)\in\mathcal{Q}^{r_l}_1$ then the event
$\mathcal{B}_i(Q_l)$ is that the subgraph induced on the $i$-th
collection is a strict $(Q_l,\{x\})$-extension of the first vertex
of the collection and forms with it a $j$-maximal pair. By
Lemma~\ref{balance-uho} for any $Q_i\in\mathcal{Q}^{r_i}_{l_i},$
$i\in\{1,...,t\}$, there exist numbers $q(Q_1,...,Q_t)>0$ such
that
$$
 \sum_{r_1+...+r_t=n}\sum_{i=1}^t\sum_{(Q_i,x)\in\mathcal{Q}^{r_i}_1}
 q(Q_1,...,Q_t)\sum_{i_1,...,i_t}{\sf P}_{N,p}
 (\mathcal{B}_{i_1}(Q_1)\cap...\cap\mathcal{B}_{i_t}(Q_t))\leq
 r_j(n,N)\leq
$$
\begin{equation}
 \leq\sum_{r_1+...+r_t=n}\sum_{i=1}^t\sum_{l_i=1}^2\sum_{(Q_i,x)\in\mathcal{Q}^{r_i}_{l_i}}
 q(Q_1,...,Q_t)\sum_{i_1,...,i_t}{\sf P}_{N,p}
 (\mathcal{B}_{i_1}(Q_1)\cap...\cap\mathcal{B}_{i_t}(Q_t)).
\label{bounds}
\end{equation}

Summations in \ref{bounds} are over $i_1,...,i_t$ corresponding to
pairwise disjoint collections.

Let $r_1,...,r_t\in\mathbb{N}$, $r_1+...+r_t=n$. Consider a vector
$(l_1,...,l_t)\in\{1,2\}^t$ such that al least one of the numbers
$l_1,...,l_t$ equals $2$. Let $(Q_i,x)\in\mathcal{Q}^{r_i}_{l_i}$.
Consider graphs $\widehat{Q}_1,...,\widehat{Q}_t$ with the
following properties.

\begin{itemize}

\item[---] Any two graphs among $\widehat{Q}_1,...,\widehat{Q}_t$
do not have a common vertex.

\item[---] There exist vertices $x_1,...,x_t$ such that
$(\widehat{Q}_i,x_i)\cong(Q_i,x)$ for any $i\in\{1,...,t\}$.

\end{itemize}

Set $V(G)=V(\widehat{Q}_1)\cup...\cup V(\widehat{Q}_t)$,
$E(G)=E(\widehat{Q}_1)\cup...\cup E(\widehat{Q}_t)$. Let $N_G$ be
the number of copies of $G$ in $G(N,p)$. Obviously there exist
$C(G),\mu(G)>0$ such that
$$
{\sf E}_{N,p} N_G<C(G) N^{v(G)-\alpha\cdot e(G)}<C(G) N^{-\mu(G)}.
$$

Therefore the difference between the upper and the lower bound in
(\ref{bounds}) equals $o(1)$. The existence of
$\lim\limits_{N\rightarrow\infty}r_j(n,N)$ follows from the
convergence of the lower bound in (\ref{bounds}) to a number
$r_j(n)$ as $N\rightarrow\infty$. Let us prove this convergence.
Let $R_1,...,R_t$ be pairwise disjoint subsets of $V_N$,
$|R_i|=|V(Q_i)|$. Let $\varphi(Q_1,...,Q_t)$ be the probability
that the graphs induced on $R_1,...,R_t$ form with the first
vertices of $R_1,...,R_t$ $j$-maximal pairs under the following
condition. For each $l\in\{1,...,t\}$ the graph induced on the
$l$-th collection is a strict $(Q_l,x)$-extension of the first
vertex of this collection. The proof of the convergence of
$\varphi(Q_1,...,Q_t)$, $N\rightarrow\infty$, is identical to the
proof of Statement~\ref{probab_converge}. Thus we do not give this
proof. The existence of $r_j(n)$ follows from the convergence of
$\varphi(Q_1,...,Q_t)$.

Finally let us prove that $r_j(n)\rightarrow 0$ when
$n\rightarrow\infty$. It is easy to see that
$$
 \sum_{r_1+...+r_t=n}\sum_{i=1}^t\sum_{(Q_i,x)\in\mathcal{Q}^{r_i}_1}
 q(Q_1,...,Q_t)\sum_{i_1,...,i_t}{\sf P}_{N,p}
 (\mathcal{B}_{i_1}(Q_1)\cap...\cap\mathcal{B}_{i_t}(Q_t))
$$
$$
 \leq
 (C_N^{v_j})^n \left(\frac{v_j!}{a_j}\right)^n\frac{1}{n!} p^{e(Q_1)+...+e(Q_t)}
 \varphi(Q_1,...,Q_t)
 \leq \frac{1}{a_j^n n!}=o(1).
$$
Therefore the convergence of ${\sf P}_{N,p}(\mathcal{L}_j^k(N))$ is proved.\\

Let us consider intersections of the properties
$\mathcal{L}_j^k(N)$.

The convergence~(\ref{A}) follows from the existence of a limit of
the sequence $\{{\sf P}_{N,p}(\mathcal{L}^k_{j_1}(N)\cap...\cap
\mathcal{L}^k_{j_t}(N))\}_{N\in\mathbb{N}}$ for any
$j_1,...,j_t\in\{1,...,m(k)\}$. Indeed, for any properties $A,C$
the equality ${\sf P}(\overline{A}\cap C)={\sf P}(C)-{\sf P}(A\cap
C)$ holds. If ${\sf
P}(\overline{A_1}\cap...\cap\overline{A_{k-1}}\cap C)$ equals
$$
 \sum_s\sum_{i_1,...,i_s}(-1)^{\sigma_s(i_1,...,i_s)}{\sf P}
 (A_{i_1}\cap...\cap A_{i_s}\cap C)
$$
for some $\sigma_s:\mathbb{N}^s\rightarrow\{0,1\}$, then
$$
 {\sf P}(\overline{A_1}\cap...\cap\overline{A_k}\cap C)=
 {\sf P} (\overline{A_1}\cap...\cap\overline{A_{k-1}}\cap
 C)-{\sf P}(\overline{A_1}\cap...\cap(A_k\cap C))=
$$
$$
 \sum_s\sum_{i_1,...,i_s}(-1)^{\sigma(i_1,...,i_s)}({\sf P}
 (A_{i_1}\cap...\cap A_{i_s}\cap C)-{\sf P}
 (A_{i_1}\cap...\cap A_{i_s}\cap A_k\cap C)).
$$
In other words, the probability ${\sf
P}(\overline{A_1}\cap...\cap\overline{A_k}\cap C)$ can be written
as the finite sum of the probabilities of some intersections of
properties without any negations. Therefore the existence of a
limit of any such intersection implies the existence of a limit of
${\sf
P}(\overline{A_1}\cap...\cap\overline{A_k}\cap C)$.\\

Thus we have $j_1,...,j_t\in\{1,...,m(k)\}$. The proof of the
existence of $\lim_{n\rightarrow\infty}{\sf
P}_{N,p}(\mathcal{L}^k_{j_1}(N)\cap...\cap\mathcal{L}^k_{j_t}(N))$
and the proof of the convergence of the probability of one
property are the same. Note that if an intersection of
$\{j_1,...,j_t\}$ and $\{m(k)+1,...,\widehat{m}(k)\}$ is not
empty, then the probability of the existence of
$X_{j_i}^{**}(\widetilde{x}_i)$ converges due to arguments which
are the same as in the case $j_1,...,j_t\in\{1,...,m(k)\}$.
Therefore it remains to apply Theorem 4. Finally the convergence
${\sf P}_{N,p}(\mathcal{L}^k_{\widehat{m}(k)+1}(N))$ follows from
the equality
$$
 \mathcal{L}^k_{\widehat{m}(k)+1}(N)=\Omega_N\setminus\mathcal{L}^k_1(N)\cup...\cup\mathcal{L}^k_{\widehat{m}(k)}(N).
$$
Lemma is proved .\\

\subsection{ Proof of Lemma~\ref{classes}}

Let $\mathcal{S}$ be the set of all $\alpha$-rigid and
$\alpha$-neutral pairs $(K,T)$ such that $v(T)\leq k^3, v(K,T)\leq
k^3$.
Theorem~\ref{maximal_extensions_rigid},
Theorem~\ref{maximal_extensions_neutral}, Lemma~\ref{max_graph}
imply the existence of a set $\widetilde{\Omega}_N\subset\Omega_N$
such that
$$
 \lim_{N\rightarrow\infty}{\sf P}_{N,p}(\Omega_N
 \setminus\widetilde{\Omega}_N)=0
$$
and the following property holds. For any
$\mathcal{G}\in\widetilde{\Omega}_N,$ $r\leq k$,
$\widetilde{x}_1,...,\widetilde{x}_r$, $(K,T)\in\mathcal{S}$ there
exist all possible non-isomorphic $(K,T)$-maximal $\alpha$-safe
pairs
$(W,\mathcal{G}|_{\{\widetilde{x}_1,...,\widetilde{x}_r\}})$,
$v(W)\leq k^3,$ in $\mathcal{G}$ and all possible non-isomorphic
$(K,T)$-maximal in $\mathcal{G}$ strictly balanced graphs $W$ with
$\rho(W)<\alpha$, $v(W)\leq k^3$, and there is no copy of a graph
with $r\leq
k^3$ vertices and density greater than $\alpha$.\\

If in some rounds a strategy of Spoiler doesn't depend on the
choice between the graphs $\mathcal{G}$ and $\mathcal{H}$ then we
assume that he chooses the graph $\mathcal{G}$.

Let us prove that for any $N,M\in\mathbb{N}$ and any pair
$$
 (\mathcal{G},\mathcal{H})\in
 (\mathcal{A}_{j_1...j_t}(N)\cap\widetilde{\Omega}_N)\times
 (\mathcal{A}_{j_1...j_t}(M)\cap \widetilde{\Omega}_M)
$$
Duplicator has a winning strategy in the game
$EHR(\mathcal{G},\mathcal{H},k)$.\\

Let Spoiler choose a vertex $\widetilde{x}_1$ in $\mathcal{G}$ at
the first round. Consider the graph
$X^5_{\mathcal{G}}(\widetilde{x}_1)$. If
$\rho^{\max}(X^5_{\mathcal{G}}(\widetilde{x}_1))=k-2$ and
$X^5_{\mathcal{G}}(\widetilde{x}_1)$ is $(\widetilde{x}_1)$-net in
$\mathcal{G}$,  then as $\mathcal{H}\in\mathcal{A}_{j_1...j_t}(M)$
there is a vertex $\widetilde{y}_1$ in $\mathcal{H}$ such that
$X^5_{\mathcal{H}}(\widetilde{y}_1)$ is a net of the graph
$X^5_{\mathcal{G}}(\widetilde{x}_1)$ (in the considered case
graphs $X^5_{\mathcal{G}}(\widetilde{x}_1)$ and
$X^5_{\mathcal{H}}(\widetilde{y}_1)$ are isomorphic). Let either
$\rho^{\max}(X^5_{\mathcal{G}}(\widetilde{x}_1))=k-2$ and
$X^5_{\mathcal{G}}(\widetilde{x}_1)$ be not a
$(\widetilde{x}_1)$-net in $\mathcal{G}$ or
$\rho^{\max}(X^5_{\mathcal{G}}(\widetilde{x}_1))<k-2$. 
As $\mathcal{H}\in\widetilde{\Omega}_M$ in the graph $\mathcal{H}$
there is a vertex $\widetilde{y}_1$ such that the graph
$X^5_{\mathcal{H}}(\widetilde{y}_1)$ is a net of the graph
$X^5_{\mathcal{G}}(\widetilde{x}_1)$. Duplicator chooses the
vertex $\widetilde{y}_1$ at the first round.

Let at the $\xi$-th round, $\xi\in\{2,...,k-3\}$, Spoiler choose a
vertex $\widetilde{x}_{\xi}\in\mathcal{G}$. If for some
$i\in\{1,t(\widetilde{x}_1,...,\widetilde{x}_{\xi-1},\mathcal{G})\}$,
$\mu\in\{\xi,...,k-3\}$ vertices $\widetilde{x}_{\xi}$ and
$\widetilde{x}^{i}_{\mu}(\widetilde{x}_1,...,\widetilde{x}_{\xi-1})$
coincide then Duplicator chooses the vertex
$$
 \widetilde{y}_{\xi}=\mbox{NET}_{\mathcal{G},\mathcal{H},
 \widetilde{x}_1,...,\widetilde{x}_{\xi-1},\widetilde{y}_1,...,\widetilde{y}_{\xi-1}}(\widetilde{x}^{i}_{\mu}
 (\widetilde{x}_1,...,\widetilde{x}_{\xi-1})).
$$

Suppose there are no appropriate
$i\in\{1,t(\widetilde{x}_1,...,\widetilde{x}_{\xi-1},\mathcal{G})\}$
and $\mu\in\{\xi,...,k-3\}$. As $\mathcal{H}\in\Omega_M$ from the
definitions of
$X^5_{\mathcal{G}}(\widetilde{x}_1,...,\widetilde{x}_{\xi-1}),$
$X^5_{\mathcal{G}}(\widetilde{x}_1,...,\widetilde{x}_{\xi})$ it
follows that in the graph $\mathcal{H}$ there is a vertex
$\widetilde{y}_{\xi}$ such that the graph
$X^5_{\mathcal{H}}(\widetilde{y}_1,...,\widetilde{y}_{\xi})$ is a
net of the graph
$X^5_{\mathcal{G}}(\widetilde{x}_1,...,\widetilde{x}_{\xi})$.
Indeed, we want to construct a graph
$X^5_{\mathcal{H}}(\widetilde{y}_1,...,\widetilde{y}_{\xi})$ such
that the pair
$(X^5_{\mathcal{H}}(\widetilde{y}_1,...,\widetilde{y}_{\xi},\mathcal{H}|_{\{\widetilde{y}_1,...,\widetilde{y}_{\xi-1}\}}))$
is $\alpha$-safe. Duplicator chooses the vertex
$\widetilde{y}_{\xi}$.

Let at the $k-3$-th round vertices
$\widetilde{x}_1,...,\widetilde{x}_{k-3}\in\mathcal{G}$,
$\widetilde{y}_1,...,\widetilde{y}_{k-3}\in\mathcal{H}$ be chosen.
The graphs
$X^5_{\mathcal{G}}(\widetilde{x}_1,...,\widetilde{x}_{k-3}),
X^5_{\mathcal{H}}(\widetilde{y}_1,...,\widetilde{y}_{k-3})$ are
isomorphic (it follows from the ways of their constructions). Let
$\varphi:
X^5_{\mathcal{G}}(\widetilde{x}_1,...,\widetilde{x}_{k-3})\rightarrow
X^5_{\mathcal{H}}(\widetilde{y}_1,...,\widetilde{y}_{k-3})$ be an
isomorphism.

The remaining part of the proof is divided into cases. There are
some basic cases such that other cases are similar to them. Thus
we give proofs for the basic cases only. However we give brief
proofs for all the cases.

\begin{itemize}

\item[1.] At the $k-2$-th round Spoiler chooses a vertex
$\widetilde{x}_{k-2}$ adjacent to
$\widetilde{x}_1,...,\widetilde{x}_{k-3}$.

If in $\mathcal{G}$ there are vertices
$\widetilde{x}^1,\widetilde{x}^2$ adjacent to each of the vertices
$\widetilde{x}_1,...,\widetilde{x}_{k-2}$, then the vertex
$\widetilde{x}_{k-2}$ is in
$V(\widehat{X}_{\mathcal{G}}^2(\widetilde{x}_1,...,\widetilde{x}_{k-3}))$.
Duplicator chooses
$\widetilde{y}_{k-2}=\varphi(\widetilde{x}_{k-2})$. If at the
$k-1$-th round Spoiler chooses a vertex from
$V(\widehat{X}_{\mathcal{G}}^2(\widetilde{x}_1,...,\widetilde{x}_{k-3}))$,
then Duplicator chooses the vertex
$\widetilde{y}_{k-1}=\varphi(\widetilde{x}_{k-1})$ again and
obviously wins. If Spoiler chooses a vertex $\widetilde{x}_{k-1}$
adjacent to each of the vertices
$\widetilde{x}_1,...,\widetilde{x}_{k-2}$ and there is no vertex
$\widetilde{x}$ adjacent to
$\widetilde{x}_1,...,\widetilde{x}_{k-1}$, then without loss of
generality one can consider $\widetilde{x}_{k-1}$ to be an element
of the set
$X^2_{\mathcal{G}}(\widetilde{x}_1,...,\widetilde{x}_{k-3})$.
Spoiler chooses the vertex
$\widetilde{y}_{k-1}=\varphi(\widetilde{x}_{k-1})$. If at the
$k$-th round Spoiler chooses a vertex $\widetilde{x}_k$ adjacent
to $k-2$ vertices from $\widetilde{x}_1,...,\widetilde{x}_{k-1}$
(say, to vertices $\widetilde{x}_1,...,\widetilde{x}_{k-2}$), then
in $X^2_{\mathcal{G}}(\widetilde{x}_1,...,\widetilde{x}_{k-3})$
there is a vertex $\widetilde{x}$, adjacent to
$\widetilde{x}_1,...,\widetilde{x}_{k-2}$. Duplicator chooses
$\widetilde{y}_k=\varphi(\widetilde{x}_k)$ and wins. Finally, if
the vertex $\widetilde{x}_k$ is adjacent to at most $k-3$ vertices
from $\widetilde{x}_1,...,\widetilde{x}_{k-1}$, then the pair
$(\mathcal{G}|_{\{\widetilde{x}_1,...,\widetilde{x}_k\}},\mathcal{G}|_{\{\widetilde{x}_1,...,\widetilde{x}_{k-1}\}})$
is $\alpha$-safe. As $\mathcal{H}\in\widetilde{\Omega}_M$, in
$\mathcal{H}$ there is a vertex $\widetilde{y}_k$ such that the
graph $\mathcal{H}|_{\{\widetilde{y}_1,...,\widetilde{y}_k\}}$ is
a strict
$(\mathcal{G}|_{\{\widetilde{x}_1,...,\widetilde{x}_k\}},\mathcal{G}|_{\{\widetilde{x}_1,...,\widetilde{x}_{k-1}}\})$-extension
of the graph
$\mathcal{H}|_{\{\widetilde{y}_1,...,\widetilde{y}_{k-1}\}}$.
Duplicator chooses $\widetilde{y}_k$ and wins. Let the vertex
$\widetilde{x}_{k-1}$ be adjacent to at most $k-4$ vertices from
$\widetilde{x}_1,...,\widetilde{x}_{k-2}$. Let
$\widetilde{x}^1,...,\widetilde{x}^s$ be all vertices from
$\mathcal{G}$ adjacent to at most $k-2$ vertices from
$\widetilde{x}_1,...,\widetilde{x}_{k-1}$ such that pairs
$(\mathcal{G}|_{\{\widetilde{x}_1,...,\widetilde{x}_{k-1},\widetilde{x}^i\}},
\mathcal{G}|_{\{\widetilde{x}_1,...,\widetilde{x}_{k-1}\}})$,
$i\in\{1,...,s\},$ are non-isomorphic. Consider a subgraph $A$ of
$\mathcal{G}$ containing the vertices
$\widetilde{x}_1,...,\widetilde{x}_{k-1},\widetilde{x}^1,...,\widetilde{x}^s$
only. The pair
$(A,\mathcal{G}|_{\{\widetilde{x}_1,...,\widetilde{x}_{k-2}\}})$
is $\alpha$-safe.  As $\mathcal{H}\in\widetilde{\Omega}_M$, in
$\mathcal{H}$ there is a strict
$(A,\mathcal{G}|_{\{\widetilde{x}_1,...,\widetilde{x}_{k-2}\}})$-extension
$B$ of the graph
$\mathcal{H}|_{\{\widetilde{y}_1,...,\widetilde{y}_{k-3}\}}$. Let
$\xi:\,A\rightarrow B$ be an isomorphism corresponding to this
extension. Duplicator chooses $\xi(\widetilde{x}_{k-1})$ and wins.
Lastly, let the vertex $\widetilde{x}_{k-1}$ be adjacent to $k-3$
vertices from $\widetilde{x}_1,...,\widetilde{x}_{k-2}$. If there
is a vertex $\widetilde{x}$ adjacent to each of
$\widetilde{x}_1,...,\widetilde{x}_{k-1}$, then without loss of
generality one can consider vertices
$\widetilde{x}_{k-1},\widetilde{x}$ to be in
$V(X^2_{\mathcal{G}}(\widetilde{x}_1,...,\widetilde{x}_{k-3}))$.
Duplicator chooses a vertex
$\widetilde{y}_{k-1}=\varphi(\widetilde{x}_{k-1})$ and wins. If
Spoiler chooses a vertex adjacent to each of
$\widetilde{x}_1,...,\widetilde{x}_{k-1}$, then Duplicator chooses
$\varphi(\widetilde{x})$. Let Spoiler choose a vertex adjacent to
$k-2$ vertices from $\widetilde{x}_1,...,\widetilde{x}_{k-1}$. We
can assume that this vertex is from
$V(X^2_{\mathcal{G}}(\widetilde{x}_1,...,\widetilde{x}_{k-3}))$.
Therefore, Duplicator wins again.

Let $\widetilde{x}_{k-2}$ be not from
$V(\widehat{X}^2_{\mathcal{G}}(\widetilde{x}_1,...,\widetilde{x}_{k-3}))$.
Let $\widetilde{x}^1,...,\widetilde{x}^s$ be all vertices of the
graph $\mathcal{G}$ such that the following properties hold. Each
of the vertices $\widetilde{x}_1,...,\widetilde{x}_{k-2}$ is
adjacent to each of the vertices
$\widetilde{x}^1,...,\widetilde{x}^s$. Sets of collections
containing $k-3$ vertices from
$\widetilde{x}_1,...,\widetilde{x}_{k-2}$ such that in
$\mathcal{G}$ there is a vertex adjacent to them and to
$\widetilde{x}^i$ are different for all  $i\in\{1,...,s\}$.
Consider a subgraph $A$ of $\mathcal{G}$ containing the vertices
$\widetilde{x}_1,...,\widetilde{x}_{k-2},\widetilde{x}^1,...,\widetilde{x}^s$
and one $(G_2,H_2)$-extension for each subgraph with $k-3$
vertices from $\widetilde{x}_1,...,\widetilde{x}_{k-2}$ and one
from $\widetilde{x}^1,...,\widetilde{x}^s$. Then the pair
$(A,\mathcal{G}|_{\{\widetilde{x}_1,...,\widetilde{x}_{k-3}\}})$
is $\alpha$-safe.  As $\mathcal{H}\in\widetilde{\Omega}_M$ in
$\mathcal{H}$ there is a strict
$(A,\mathcal{G}|_{\{\widetilde{x}_1,...,\widetilde{x}_{k-3}\}})$-extension
$B$ of the graph
$\mathcal{H}|_{\{\widetilde{y}_1,...,\widetilde{y}_{k-3}\}}$. Let
$\xi:\,A\rightarrow B$ be an isomorphism corresponding to this
extension. Duplicator chooses
$\widetilde{y}_{k-2}=\xi(\widetilde{x}_{k-2})$. Let at the
$k-1$-th round Spoiler choose a vertex $\widetilde{x}_{k-1}$
adjacent to each of $\widetilde{x}_1,...,\widetilde{x}_{k-2}$.
There is such a vertex in $B$ that Duplicator can win by choosing
this vertex. If Spoiler chooses a vertex $\widetilde{x}_{k-1}$
adjacent to at most $k-4$ vertices from
$\widetilde{x}_1,...,\widetilde{x}_{k-2}$, then obviously
Duplicator has a winning strategy. Finally, let
$\widetilde{x}_{k-1}$ be adjacent to $k-3$ vertices from
$\widetilde{x}_1,...,\widetilde{x}_{k-2}$. If there is a vertex
$\widetilde{x}$ adjacent to each of
$\widetilde{x}_1,...,\widetilde{x}_{k-1}$, then without loss of
generality we can consider vertices
$\widetilde{x}_{k-1},\widetilde{x}$ to be in $V(A)$. Then there is
a vertex $\widetilde{y}_{k-1}$ in $B$, corresponding to the vertex
$\widetilde{x}_{k-1}$. Duplicator chooses this vertex and wins.

\item[2.] At the $k-2$-th round Spoiler chooses a vertex
$\widetilde{x}_{k-2}$ adjacent to $k-4$ vertices from
$\widetilde{x}_1,...,\widetilde{x}_{k-3}$.

If in $\mathcal{G}$ there are vertices
$\widetilde{x}^1,\widetilde{x}^2,$ adjacent to each of
$\widetilde{x}_1,...,\widetilde{x}_{k-2}$, then either the
vertices $\widetilde{x}^1,\widetilde{x}^2$ are in
$V(\widehat{X}_{\mathcal{G}}^2(\widetilde{x}_1,...,\widetilde{x}_{k-3}))$,
or the vertices $\widetilde{x}^1,\widetilde{x}^2$ are in
$V(\widehat{X}_{\mathcal{G}}^4(\widetilde{x}_1,...,\widetilde{x}_{k-3}))\setminus
V(\widehat{X}_{\mathcal{G}}^2(\widetilde{x}_1,...,\widetilde{x}_{k-3}))$.
Anyway the vertex $\widetilde{x}_{k-2}$ is in
$V(X_{\mathcal{G}}^4(\widetilde{x}_1,...,\widetilde{x}_{k-3}))$.
Spoiler chooses
$\widetilde{y}_{k-2}=\varphi(\widetilde{x}_{k-2})$. Further
choices of Duplicator are described in the same manner as for the
case 1.

If there are no two vertices adjacent to each other and to each of
$\widetilde{x}_1,...,\widetilde{x}_{k-2}$, then further reasonings
are identical to reasonings from subcases of the case 1., in which
we use safe pairs.

\item[3.] At the $k-2$-th round Spoiler chooses a vertex
$\widetilde{x}_{k-2}$, adjacent to at most $k-5$ vertices from
$\widetilde{x}_1,...,\widetilde{x}_{k-3}$.

If in $\mathcal{G}$ there are vertices
$\widetilde{x}^1,\widetilde{x}^2,$ adjacent to each of
$\widetilde{x}_1,...,\widetilde{x}_{k-2}$, then the vertices
$\widetilde{x}^1,\widetilde{x}^2$ are in one of the sets
$V(\widehat{X}_{\mathcal{G}}^2(\widetilde{x}_1,...,\widetilde{x}_{k-3}))$,
$V(\widehat{X}_{\mathcal{G}}^4(\widetilde{x}_1,...,\widetilde{x}_{k-3}))\setminus
V(\widehat{X}_{\mathcal{G}}^2(\widetilde{x}_1,...,\widetilde{x}_{k-3}))$,
$V(\widehat{X}_{\mathcal{G}}^5(\widetilde{x}_1,...,\widetilde{x}_{k-3}))\setminus
V(\widehat{X}_{\mathcal{G}}^4(\widetilde{x}_1,...,\widetilde{x}_{k-3}))$.
Anyway without loss of generality we can consider the vertex
$\widetilde{x}_{k-2}$ to be in
$V(X_{\mathcal{G}}^5(\widetilde{x}_1,...,\widetilde{x}_{k-3}))$.
Spoiler chooses
$\widetilde{y}_{k-2}=\varphi(\widetilde{x}_{k-2})$. Further
choices of Duplicator are described in the same manner as for the
case 1.

If there are no two vertices adjacent to each other and to each of
$\widetilde{x}_1,...,\widetilde{x}_{k-2}$, then further reasonings
are identical to reasonings from subcases of the case 1., in which
we use safe pairs.

\end{itemize}

Lemma is proved.

\end{document}